\documentclass{article}
\usepackage{verbatim,amsmath,amscd,amssymb,amsthm}
\usepackage[all]{xy}
\begin{document}
\font\germ=eufm10
\def\ssl{\hbox{\germ sl}}
\def\slh{\widehat{\ssl_2}}
\def\ge{\hbox{\germ g}}

\makeatletter
\def\aaa{@}
\centerline{}
\centerline{\Large\bf Geometric Crystals 
on Schubert Varieties}
\vskip7pt
\vskip15pt
\centerline{NAKASHIMA Toshiki}
\vskip10pt
\centerline{Department of Mathematics,}
\centerline{Sophia University, Tokyo 102-8554, JAPAN}
\centerline{e-mail:\,\,toshiki@mm.sophia.ac.jp}
\vskip10pt
\centerline{Abstract}
We define geometric crystals and 
unipotent crystals for 
arbitrary Kac-Moody groups and 
describe geometric and unipotent 
crystal structures on the Schubert varieties.

\makeatother

\renewcommand{\labelenumi}{$($\roman{enumi}$)$}
\renewcommand{\labelenumii}{$(${\rm \alph{enumii}}$)$}
\font\germ=eufm10
\newcommand{\cI}{{\mathcal I}}
\newcommand{\cA}{{\mathcal A}}
\newcommand{\cB}{{\mathcal B}}
\newcommand{\cC}{{\mathcal C}}
\newcommand{\cD}{{\mathcal D}}
\newcommand{\cF}{{\mathcal F}}
\newcommand{\cH}{{\mathcal H}}
\newcommand{\cK}{{\mathcal K}}
\newcommand{\cL}{{\mathcal L}}
\newcommand{\cM}{{\mathcal M}}
\newcommand{\cN}{{\mathcal N}}
\newcommand{\cO}{{\mathcal O}}
\newcommand{\cS}{{\mathcal S}}
\newcommand{\cV}{{\mathcal V}}
\newcommand{\fra}{\mathfrak a}
\newcommand{\frb}{\mathfrak b}
\newcommand{\frc}{\mathfrak c}
\newcommand{\frd}{\mathfrak d}
\newcommand{\fre}{\mathfrak e}
\newcommand{\frf}{\mathfrak f}
\newcommand{\frg}{\mathfrak g}
\newcommand{\frh}{\mathfrak h}
\newcommand{\fri}{\mathfrak i}
\newcommand{\frj}{\mathfrak j}
\newcommand{\frk}{\mathfrak k}
\newcommand{\frI}{\mathfrak I}
\newcommand{\fm}{\mathfrak m}
\newcommand{\frn}{\mathfrak n}
\newcommand{\frp}{\mathfrak p}
\newcommand{\fq}{\mathfrak q}
\newcommand{\frr}{\mathfrak r}
\newcommand{\frs}{\mathfrak s}
\newcommand{\frt}{\mathfrak t}
\newcommand{\fru}{\mathfrak u}
\newcommand{\frA}{\mathfrak A}
\newcommand{\frB}{\mathfrak B}
\newcommand{\frF}{\mathfrak F}
\newcommand{\frG}{\mathfrak G}
\newcommand{\frH}{\mathfrak H}
\newcommand{\frJ}{\mathfrak J}
\newcommand{\frN}{\mathfrak N}
\newcommand{\frP}{\mathfrak P}
\newcommand{\frT}{\mathfrak T}
\newcommand{\frU}{\mathfrak U}
\newcommand{\frV}{\mathfrak V}
\newcommand{\frX}{\mathfrak X}
\newcommand{\frY}{\mathfrak Y}
\newcommand{\frZ}{\mathfrak Z}
\newcommand{\rA}{\mathrm{A}}
\newcommand{\rC}{\mathrm{C}}
\newcommand{\rd}{\mathrm{d}}
\newcommand{\rB}{\mathrm{B}}
\newcommand{\rD}{\mathrm{D}}
\newcommand{\rE}{\mathrm{E}}
\newcommand{\rH}{\mathrm{H}}
\newcommand{\rK}{\mathrm{K}}
\newcommand{\rL}{\mathrm{L}}
\newcommand{\rM}{\mathrm{M}}
\newcommand{\rN}{\mathrm{N}}
\newcommand{\rR}{\mathrm{R}}
\newcommand{\rT}{\mathrm{T}}
\newcommand{\rZ}{\mathrm{Z}}
\newcommand{\bbA}{\mathbb A}
\newcommand{\bbC}{\mathbb C}
\newcommand{\bbG}{\mathbb G}
\newcommand{\bbF}{\mathbb F}
\newcommand{\bbH}{\mathbb H}
\newcommand{\bbP}{\mathbb P}
\newcommand{\bbN}{\mathbb N}
\newcommand{\bbQ}{\mathbb Q}
\newcommand{\bbR}{\mathbb R}
\newcommand{\bbV}{\mathbb V}
\newcommand{\bbZ}{\mathbb Z}
\newcommand{\adj}{\operatorname{adj}}
\newcommand{\Ad}{\mathrm{Ad}}
\newcommand{\Ann}{\mathrm{Ann}}
\newcommand{\rcris}{\mathrm{cris}}
\newcommand{\ch}{\mathrm{ch}}
\newcommand{\coker}{\mathrm{coker}}
\newcommand{\diag}{\mathrm{diag}}
\newcommand{\Diff}{\mathrm{Diff}}
\newcommand{\Dist}{\mathrm{Dist}}
\newcommand{\rDR}{\mathrm{DR}}
\newcommand{\ev}{\mathrm{ev}}
\newcommand{\Ext}{\mathrm{Ext}}
\newcommand{\cExt}{\mathcal{E}xt}
\newcommand{\fin}{\mathrm{fin}}
\newcommand{\Frac}{\mathrm{Frac}}
\newcommand{\GL}{\mathrm{GL}}
\newcommand{\Hom}{\mathrm{Hom}}
\newcommand{\hd}{\mathrm{hd}}
\newcommand{\rht}{\mathrm{ht}}
\newcommand{\id}{\mathrm{id}}
\newcommand{\im}{\mathrm{im}}
\newcommand{\inc}{\mathrm{inc}}
\newcommand{\ind}{\mathrm{ind}}
\newcommand{\coind}{\mathrm{coind}}
\newcommand{\Lie}{\mathrm{Lie}}
\newcommand{\Max}{\mathrm{Max}}
\newcommand{\mult}{\mathrm{mult}}
\newcommand{\op}{\mathrm{op}}
\newcommand{\ord}{\mathrm{ord}}
\newcommand{\pt}{\mathrm{pt}}
\newcommand{\qt}{\mathrm{qt}}
\newcommand{\rad}{\mathrm{rad}}
\newcommand{\res}{\mathrm{res}}
\newcommand{\rgt}{\mathrm{rgt}}
\newcommand{\rk}{\mathrm{rk}}
\newcommand{\SL}{\mathrm{SL}}
\newcommand{\soc}{\mathrm{soc}}
\newcommand{\Spec}{\mathrm{Spec}}
\newcommand{\St}{\mathrm{St}}
\newcommand{\supp}{\mathrm{supp}}
\newcommand{\Tor}{\mathrm{Tor}}
\newcommand{\Tr}{\mathrm{Tr}}
\newcommand{\wt}{\mathrm{wt}}
\newcommand{\Ab}{\mathbf{Ab}}
\newcommand{\Alg}{\mathbf{Alg}}
\newcommand{\Grp}{\mathbf{Grp}}
\newcommand{\Mod}{\mathbf{Mod}}
\newcommand{\Sch}{\mathbf{Sch}}\newcommand{\bfmod}{{\bf mod}}
\newcommand{\Qc}{\mathbf{Qc}}
\newcommand{\Rng}{\mathbf{Rng}}
\newcommand{\Top}{\mathbf{Top}}
\newcommand{\Var}{\mathbf{Var}}
\newcommand{\gromega}{\langle\omega\rangle}
\newcommand{\lbr}{\begin{bmatrix}}
\newcommand{\rbr}{\end{bmatrix}}
\newcommand{\for}{\bigcirc\kern-2.6ex \because}
\newcommand{\forb}{\bigcirc\kern-2.8ex \because}
\newcommand{\forbb}{\bigcirc\kern-3.0ex \because}
\newcommand{\forbbb}{\bigcirc\kern-3.1ex \because}
\newcommand{\cd}{commutative diagram }
\newcommand{\SpS}{spectral sequence}
\newcommand\C{\mathbb C}
\newcommand\hh{{\hat{H}}}
\newcommand\eh{{\hat{E}}}
\newcommand\F{\mathbb F}
\newcommand\fh{{\hat{F}}}

\def\AA{{\cal A}}
\def\al{\alpha}
\def\bq{B_q(\ge)}
\def\bqm{B_q^-(\ge)}
\def\bqz{B_q^0(\ge)}
\def\bqp{B_q^+(\ge)}
\def\beneme{\begin{enumerate}}
\def\beq{\begin{equation}}
\def\beqn{\begin{eqnarray}}
\def\beqnn{\begin{eqnarray*}}
\def\bigsl{{\hbox{\fontD \char'54}}}
\def\bbra#1,#2,#3{\left\{\begin{array}{c}\hspace{-5pt}
#1;#2\\ \hspace{-5pt}#3\end{array}\hspace{-5pt}\right\}}
\def\cd{\cdots}
\def\CC{\hbox{\bf C}}
\def\ddd{\hbox{\germ D}}
\def\del{\delta}
\def\Del{\Delta}
\def\Delr{\Delta^{(r)}}
\def\Dell{\Delta^{(l)}}
\def\Delb{\Delta^{(b)}}
\def\Deli{\Delta^{(i)}}
\def\Delre{\Delta^{\rm re}}
\def\ei{e_i}
\def\eit{\tilde{e}_i}
\def\eneme{\end{enumerate}}
\def\ep{\epsilon}
\def\eeq{\end{equation}}
\def\eeqn{\end{eqnarray}}
\def\eeqnn{\end{eqnarray*}}
\def\fit{\tilde{f}_i}
\def\FF{{\rm F}}
\def\ft{\tilde{f}}
\def\gau#1,#2{\left[\begin{array}{c}\hspace{-5pt}#1\\
\hspace{-5pt}#2\end{array}\hspace{-5pt}\right]}
\def\ge{\hbox{\germ g}}
\def\gl{\hbox{\germ gl}}
\def\hom{{\hbox{Hom}}}
\def\ify{\infty}
\def\io{\iota}
\def\kp{k^{(+)}}
\def\km{k^{(-)}}
\def\llra{\relbar\joinrel\relbar\joinrel\relbar\joinrel\rightarrow}
\def\lan{\langle}
\def\lar{\longrightarrow}
\def\max{{\rm max}}
\def\lm{\lambda}
\def\Lm{\Lambda}
\def\mapright#1{\smash{\mathop{\longrightarrow}\limits^{#1}}}
\def\mm{{\bf{\rm m}}}
\def\nd{\noindent}
\def\nn{\nonumber}
\def\nnn{\hbox{\germ n}}
\def\catob{{\cal O}(B)}
\def\oint{{\cal O}_{\rm int}(\ge)}
\def\ot{\otimes}
\def\op{\oplus}
\def\opi{\ovl\pi_{\lm}}
\def\ovl{\overline}
\def\plm{\Psi^{(\lm)}_{\io}}
\def\qq{\qquad}
\def\q{\quad}
\def\qed{\hfill\framebox[2mm]{}}
\def\QQ{\hbox{\bf Q}}
\def\qi{q_i}
\def\qii{q_i^{-1}}
\def\ra{\rightarrow}
\def\ran{\rangle}
\def\rlm{r_{\lm}}
\def\ssl{\hbox{\germ sl}}
\def\slh{\widehat{\ssl_2}}
\def\ti{t_i}
\def\tii{t_i^{-1}}
\def\til{\tilde}
\def\tm{\times}
\def\tt{{\hbox{\germ{t}}}}
\def\ttt{\hbox{\germ t}}
\def\ua{U_{\AA}}
\def\ue{U_{\vep}}
\def\uq{U_q(\ge)}
\def\ufin{U^{\rm fin}_{\vep}}
\def\ufinp{(U^{\rm fin}_{\vep})^+}
\def\ufinm{(U^{\rm fin}_{\vep})^-}
\def\ufinz{(U^{\rm fin}_{\vep})^0}
\def\uqm{U^-_q(\ge)}
\def\uqp{U^+_q(\ge)}
\def\uqmq{{U^-_q(\ge)}_{\bf Q}}
\def\uqpm{U^{\pm}_q(\ge)}
\def\uqq{U_{\bf Q}^-(\ge)}
\def\uqz{U^-_{\bf Z}(\ge)}
\def\ures{U^{\rm res}_{\AA}}
\def\urese{U^{\rm res}_{\vep}}
\def\uresez{U^{\rm res}_{\vep,\ZZ}}
\def\util{\widetilde\uq}
\def\uup{U^{\geq}}
\def\ulow{U^{\leq}}
\def\bup{B^{\geq}}
\def\blow{\ovl B^{\leq}}
\def\vep{\varepsilon}
\def\vp{\varphi}
\def\vpi{\varphi^{-1}}
\def\VV{{\cal V}}
\def\xii{\xi^{(i)}}
\def\Xiioi{\Xi_{\io}^{(i)}}
\def\WW{{\cal W}}
\def\wtil{\widetilde}
\def\what{\widehat}
\def\wpi{\widehat\pi_{\lm}}
\def\ZZ{\mathbb Z}

\renewcommand{\thesection}{\arabic{section}}
\section{Introduction}
\setcounter{equation}{0}
\renewcommand{\theequation}{\thesection.\arabic{equation}}

The theory of crystal base introduced by Kashiwara 
succeeds in being
 applied to many areas in mathematics and 
mathematical physics to clarify thier 
combinatorial behavior. 
One of the reasons why it
can be well-applied is that it contains
not only ``real crystals'' 
but also ``virtual crystals'',
{\it e.g.,} $B_i$ (see Sect.6), 
$t_\lm$ (see \cite{K3}),etc, 
where ``virtual crystals'' mean 
certain purly combinatorial
objects not having the corresponding $\uq$-modules.
They are deeply related to each other by 
some `limit' operations and then 
they describe many combinatorial phenomena together.
Indeed, some ``real crystal'' is obtained as 
subcrystal in infinitely many 
tensor products of ``virtual crystals''
(see \cite{N} and the references therein) 
and the crystal $B_i$ and $\til B$ in Sect.5 
are also obtained as some limit of real crystals.
In this sense, the theory of crystals would be 
wider than usual represenatation theory of the
quantum algebra $\uq$.
Roughly, we can say that real crystal bases are 
obtained by
taking the limit $q\rightarrow 0$ from 
some bases of $\uq$-modules,
which is called ``crystallization''.
But, ``virtual'' ones are not gotten 
by such crystallizations from $\uq$-modules.

Berenstein and Kazhdan clarify (\cite{BK})
that such ``virtual crystals'' also have
some ``real'' backgrounds as the 
``tropicalization/ultra-discritization''
of ``geometric crystals''
for semi-simple(reductive) groups.

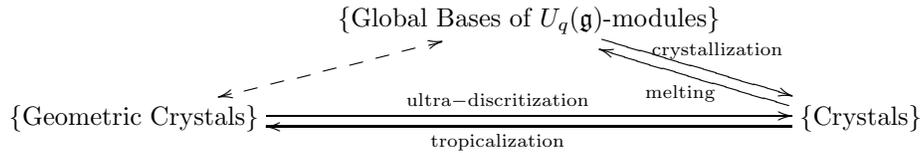
\[
\xymatrix@R=3pt{
&\{\hbox{Global Bases of $\uq$-modules}\}
\ar@{<-->}[dddl]
\ar@{->}^{
\rm\scriptstyle crystallization} @<1pt>[dddr]
\ar@{<-}_{
\rm\scriptstyle melting} @<-3pt>[dddr]
& \\ \\ \\
\{\rm Geometric\,\, Crystals\}\ar@{->}^{
\rm\scriptstyle ultra-discritization} @<1pt>[rr]
\ar@{<-}_{
\rm\scriptstyle tropicalization} @<-3pt>[rr]
&& \{\rm Crystals\} \\
}
\]

Recently, by the ultra-discritization/tropicalization
 method, the relations between soliton cellular
automaton and crystals are revealed 
(see {\it e.g.,}\cite{HH}\cite{HK}).
In the meanwhile, it is well-known that 
flag varieties $G/B$(reps.  $G/P$)
plays a significant role in the soliton theory,
where $G$ is an affine Kac-Moody group
and $B$ (resp. $P$)
is its Borel (resp. parabolic) subgroup.
We would like to find the connection of 
affine flag varieties and  geometric crystals.
For the purpose, we shall extend the
theory of geometric/unipotent 
 crystals \cite{BK} to Kac-Moody setting.
And then we shall 
define geomtric/unipotent crystals
on finite Schubert cells/varieties and 
consider some positive structures on them.
Finally, we show that some ultra-discritizations
of the geometric crystals on Schubert varieties are
isomorphic to tensor products of some Kashiwara's crystals.
But in general, it is still unclear to define 
geometric crystal structures on flag varieties 
as mentioned in the remark in Sect.4.

The organization of the article is as follows;
in Sect.2 we review briefly the theory of Kac-Moody groups,
ind-varieties and ind-groups.
In Sect.3, we define the notion of 
unipotent crystals in Kac-Moody setting
and their product structures.
We also define the notion of geometric crystals 
and give a recipe
for obtaining canonically geometric crystals 
from unipotent crystals
following to \cite{BK}.
In Sect.4, on finite Schubert cells/varieties
we induce the structure of 
unipotent/geomtric crystals.
In Sect.5, we recall the notion of positive structure 
on geometric crystals and define ultra-discritization/tropicalization 
operations. 
We also consider certain positive strucutre of 
geometric crystals on Schubert cells and show 
that its ultra-discritization is isomorphic 
to (Langlands dual of) Kashiwara's crystal $B_{i_1}\ot\cd\ot B_{i_l}$.
In the last section, we apply the result in Sect.5
 to give a new proof of braid-type 
isomorphisms (\cite{N}).

The author would acknowledge M.Kashiwara and M.Tsuzuki for valuable 
discusions and helpful advises.

\renewcommand{\thesection}{\arabic{section}}
\section{Kac-Moody groups and Ind-varieties}
\setcounter{equation}{0}
\renewcommand{\theequation}{\thesection.\arabic{equation}}

\newtheorem{pro2}{Proposition}[section]
\theoremstyle{definition}
\newtheorem{def2}[pro2]{Definition}
\theoremstyle{plain}
\newtheorem{lem2}[pro2]{Lemma}
\newtheorem{thm2}[pro2]{Theorem}
\newtheorem{ex2}[pro2]{Example}

In this section, we review on Kac-Moody groups 
following to
\cite{KP},\cite{PK},\cite{Slo}.

\subsection{Kac-Moody algebras and Kac-Moody groups}
Fix a symmetrizable generalized Cartan matix
 $A=(a_{ij})_{i,j\in I}$, where $I$ be a finite index set.
Let $(\tt,\{\al_i\}_{i\in I},\{h_i\}_{i\in I})$ 
be the associated
root data, where ${\tt}$ be the vector space 
over $\bbC$ with 
dimension $|I|+$ corank$(A)$, and 
$\{\al_i\}_{i\in I}\subset\tt^*$ and 
$\{h_i\}_{i\in I}\subset\tt$
are linearly independent indexed sets 
satisfying $\al_i(h_j)=a_{ij}$.

The Kac-Moody Lie algebra $\ge=\ge(A)$ associatd with $A$
is the Lie algbera over $\bbC$ generated by $\tt$, the 
Chevalley generators $e_i$ and $f_i$ $(i\in I)$
with the usual defining relations (\cite{KP},\cite{PK}).
There is the root space decomposition 
$\ge=\bigoplus_{\al\in \tt^*}\ge_{\al}$.
Denote the set of roots by 
$\Delta:=\{\al\in \tt^*|\al\ne0,\,\,\ge_{\al}\ne(0)\}$.
Set $Q=\sum_i\bbZ \al_i$, $Q_+=\sum_i\bbZ_{\geq0} \al_i$
and $\Delta_+:=\Delta\cap Q_+$.
An element of $\Delta_+$ is called a positive root.

Define simple reflections $s_i\in{\rm Aut}(\tt)$ $(i\in I)$ by
$s_i(h):=h-\al_i(h)h_i$, which genrate the Weyl group $W$.
We also define the action of $W$ on $\tt^*$ by
$s_i(\lm):=\lm-\al(h_i)\al_i$.
Set $\Delre:=\{w(\al_i)|w\in W,\,\,i\in I\}$, whose element 
is called a real root.

Let $\ge'$ be the derived Lie algebra of $\ge$
and $G^*$ be the free group generated by 
the free product of the additive groups 
$\ge_{\al}$ $(\al\in \Delre)$, with the canonical inclusion
$i_{\al}:\ge_{\al}\hookrightarrow G^*$.
For any integrable $\ge'$-module $(V,\pi)$, 
a homomorphism $\pi_V^*:G^*\longrightarrow {\rm Aut}_{\bbC}(V)$
is defined by $\pi_V^*(i_{\al}(e))=\exp\pi(e)$.
Set $N^*:=\cap_{V:\,{\rm integrable}}{\rm Ker}(\pi^*_V)$
and $G:=G^*/N^*$, which is called a Kac-Moody group associated 
with the Kac-Moody Lie algebra $\ge'$.
Let $\rho:G^*\rightarrow G$ be the canonical homomorphism.
For $e\in \ge_{\al}$ $(\al\in \Delre)$, define
$\exp e:=\rho(i_{\al}(e))$ and $U_{\al}:=\exp\ge_{\al}$,
which is an one-parameter subgroup of $G$.
The group $G$ is generated by $U_{\al}$ $(\al\in \Delre)$.
Let $U^{\pm}$ be the subgroups generated by $U_{\pm\al}$
($\al\in \Delre_+=\Delre\cap Q_+$), {\it i.e.,}
$U^{\pm}:=\lan U_{\pm\al}|\al\in\Del^{\rm re}_+\ran$.

For any $i\in I$, there exists a unique homomorphism;
$\phi_i:SL_2(\bbC)\rightarrow G$ such that
\[
\phi_i\left(
\left(
\begin{array}{cc}
1&t\\
0&1
\end{array}
\right)\right)=\exp t e_i,\,\,
 \phi_i\left(
\left(
\begin{array}{cc}
1&0\\
t&1
\end{array}
\right)\right)=\exp t f_i\,(t\in\bbC).
\]
Set $G_i:=\phi_i(SL_2(\bbC))$,
$T_i:=\phi_i(\{{\rm diag}(t,t^{-1})|t\in\bbC\})$ and 
$N_i:=N_{G_i}(T_i)$. Let
$T$ (resp. $N$) be the subgroup of $G$ generated by $T_i$
(resp. $N_i$), which is called a {\it maximal torus} in $G$. 
We have the isomorphism
$\phi:W\mapright{\sim}N/T$ defined by $\phi(s_i)=N_iT/T$.
An element $\ovl s_i:=x_i(-1)y_i(1)x_i(-1)$ is in 
$N_G(T)$, which is a representative of 
$s_i\in W=N_G(T)/T$. 
Define $R(w)$ for $w\in W$ by
\[
 R(w):=\{(i_1,i_2,\cd,i_l)\in I^l|w=s_{i_1}s_{i_2}\cd s_{i_l}\},
\]
where $l$ is the length of $w$.
We associate to each $w\in W$ its standard representative 
$\bar w\in N_G(T)$ by 
\[
 \bar w=\bar s_{i_1}\bar s_{i_2}\cd\bar s_{i_l},
\]
for any $(i_1,i_2,\cd,i_l)\in R(w)$.

\subsection{Ind-variety and Ind-group}
Let us recall the notion of ind-varieties and ind-groups.
(see \cite{Ku2}).
\begin{def2}
Let $k$ be an algebraically closed field. 
\begin{enumerate}
\item 
A set $X$ is an {\it ind-variety}
 over $k$ if there exists
a filtration $X_0\subset X_1\subset X_2\subset\cd$
such that 
\begin{enumerate}
\item $\displaystyle{\bigcup}_{n\geq0}X_n=X$.
\item
Each $X_n$ is a finite-dimensional variety over $k$
such that the inclusion $X_n\hookrightarrow X_{n+1}$
is a closed embedding.
\end{enumerate}
The ring of regular functions $k[X]$ is defined by 
\[
 k[X]:=\lim_{\longleftarrow\atop{n}}k[X_n].
\]
\item
A {\it Zariski topology} on an ind-variety $X$ is defined
as follows; 
a set $U\subset X$ is open if and only 
if $U\cap X_n$ is open in $X_n$ for any $n\geq0$.
\item
Let $X$ and $Y$ be two ind-varieties with filtrations
$\{X_n\}$ and $\{Y_n\}$ respectively. A map 
$f:X\rightarrow Y$ is a {\it morphism} 
if for any $n\geq0$,
there exists $m$ such that $f(X_n)\subset Y_m$
and $f_{|X_n}:X_n\rightarrow Y_m$ is a morphism.
A morphism $f:X\ra Y$ is said to be an {\it isomorphism}
if $f$ is bijective and $f^{-1}:Y\ra X$ is also 
a morphism.
\item
Let $X$ and $Y$ be two ind-varieties. A 
{\it rational morphism}
 $f:X\ra Y$ is an equivalence class
of morphisms $f_U:U\ra Y$ where $U$ is an open dense
subset of $X$, and two morphisms 
$f_U:U\ra Y$ and $f_V:V\ra Y$ are equivalnet
if they coincide on $U\cap V$.
\end{enumerate}
\end{def2}
\begin{lem2}
\label{product}
\begin{enumerate}
\item A finite dimensional variety over $k$ 
holds canonically an ind-variety structure. 
\item
If $X$ and $Y$ are ind-varieties, then
$X\tm Y$ is canonically an ind-variety by taking 
the filtration
\[
 (X\tm Y)_n:=X_n\tm Y_n.
\]
\end{enumerate}
\end{lem2}
\begin{def2}
An  ind-variety $H$ is called an
{\it ind $($algebraic$)$-group} 
if the underlying set $H$ is a group
and the maps
\[
\begin{array}{ccc}
H\tm H &\longrightarrow &H \\
(x,y)& \mapsto & xy
\end{array}\qq
\begin{array}{ccc}
H&\longrightarrow & H\\
x&\mapsto& x^{-1}
\end{array}
\]
are morphisms of ind-varieties.
\end{def2}
\begin{pro2}[\cite{Ku2}]
\label{ind-g}
\begin{enumerate}
\item
Let $G$ be a Kac-Moody group and 
$U^\pm$, $B^\pm$ be its subgroups as above.
Then $G$ is an ind-group and $U^\pm$, $B^\pm$
are its closed ind-subgroups.
\item
The multiplication maps 
\[
\begin{array}{ccc}
T\tm U &\longrightarrow &B \\
(t,u)& \mapsto & tu
\end{array}\qq
\begin{array}{ccc}
U^-\tm T&\longrightarrow & B^-\\
(v,t)&\mapsto& vt
\end{array}
\]
are isomorphisms of ind-varieties.
\end{enumerate}
\end{pro2}
\renewcommand{\thesection}{\arabic{section}}
\section{Geometric Crystals and Unipotent Crystals}
\setcounter{equation}{0}
\renewcommand{\theequation}{\thesection.\arabic{equation}}
\theoremstyle{definition}
\newtheorem{def3}{Definition}[section]
\theoremstyle{plain}
\newtheorem{pro3}[def3]{Proposition}
\newtheorem{lem3}[def3]{Lemma}
\newtheorem{thm3}[def3]{Theorem}

In this section, we define geomtric crystals and unipotent
crystals associated with Kac-Moody groups,
which is just a generalization 
of \cite{BK} to a Kac-Moody setting.

\subsection{Geometric Crystals}
Let $(a_{ij})_{i,j\in I}$ be a symmetrizable
generalized Cartan matrix and 
$G$ be the associated Kac-Moody group with the
maximal torus $T$.
An element in  $\hom(T,\bbC^\times)$
(resp. $\hom(\bbC^\times,T)$) 
is called a {\it character}
(resp. {\it co-character}) of $T$.
We define a {\it simple co-root} 
$\al_i^\vee\in \hom(\bbC^\times,T)$
$(i\in I)$ by $\al_i^\vee(t):=T_i$. 
We have a pairing $\lan \al^\vee_i,\al_j\ran=a_{ij}$.

\begin{def3}
\begin{enumerate}
\item
Let $X$ be an ind-variety over $\bbC$, 
$\gamma:X\rightarrow T$
be a rational morphism and a family of 
rational morphisms 
$e_i:\bbC^\times \times X\rightarrow X$ 
$(i\in I)$; 
\[
\begin{array}{cccc}
 e^c_i&:\bbC^\times \times X&\longrightarrow& X\\
&(c,x)&\mapsto& e^c_i(x).
\end{array}
\]
The triplet $\chi=(X,\gamma,\{e_i\}_{i\in I})$
is a {\it geometric pre-crystal} if it 
satisfies $e^1(x)=x$ and 
\begin{equation}
\label{gamma}
\gamma(e^c_i(x))=\al_i^\vee(c)\gamma(x).
\end{equation}
\item
Let $(X,\gamma_X,\{e^X_i\}_{i\in I})$ and 
$(Y,\gamma_Y,\{e^Y_i\}_{i\in I})$ be geometric pre-crystals. 
A rational morphism $f:X\ra Y$ is a 
{\it moprhism of geometric pre-crystals}
if $f$ satisfies that 
\[
f\circ e^X_i=e^Y_i\circ f, \q 
\gamma_X=\gamma_Y\circ f.
\]
In particular, if a morphism $f$ is a birational isomorphism 
of ind-vareities, it is called an 
{\it isomorphism of geometric
pre-crystals}.
\end{enumerate}
\end{def3}
Let $\chi=(X,\gamma,\{e_i\}_{i\in I})$ be a geometric pre-crystal.
For a word ${\bf i}=(i_1,i_2,\cd,i_l)\in R(w)$ 
$(w\in W)$,
set $\al^{(l)}:=\al_{i_l}$, 
$\al^{(l-1)}:=s_{i_l}(\al_{i_{l-1}})$, $\cd$,
$\al^{(1)}:=s_{i_l}\cd s_{i_2}(\al_{i_1})$.
Now for a word ${\bf i}=(i_1,i_2,\cd,i_l)\in R(w)$
we define a rational morphism $e_{\bf i}:T\times X\rightarrow X$
by
\[
(t,x)\mapsto e_{\bf i}^t(x):=e_{i_1}^{\al^{(1)}(t)}
e_{i_2}^{\al^{(2)}(t)}\cd e_{i_l}^{\al^{(l)}(t)}(x).
\]
\begin{def3}
\begin{enumerate}
\item
A geometric pre-crystal $\chi$ is called 
a {\it geometric crystal} 
if for any $w\in W$, and any ${\bf i}$, 
${\bf i'}\in R(w)$ we have
\begin{equation}					
e_{\bf i}=e_{\bf i'}.\label{ei=ei'}
\end{equation}
\item
Let $(X,\gamma_X,\{e^X_i\}_{i\in I})$ and 
$(Y,\gamma_Y,\{e^Y_i\}_{i\in I})$ be geometric crystals. 
A rational morphism $f:X\ra Y$ is called 
a {\it morphism} (resp. an
{\it isomorphism}) 
{\it of geometric crystals} if it is a morphism
(resp. an isomorphism) of geometric pre-crystals.
\end{enumerate}
\end{def3}
The following lemma is a direct result from 
\cite{BK}[Lemma 2.1] and the fact that 
the Weyl group of any 
Kac-Moody Lie algebra is a Coxeter group 
\cite{Kac}[Proposition 3.13].
\begin{lem3}
\label{Verma}
The relations $(\ref{ei=ei'})$ 
are equivalent to the following
relations:
\[
 \begin{array}{lll}
&\hspace{-20pt}e^{c_1}_{i}e^{c_2}_{j}
=e^{c_2}_{j}e^{c_1}_{i}&
{\rm if }\,\,\lan \al^\vee_i,\al_j\ran=0,\label{0}\\
&\hspace{-20pt} e^{c_1}_{i}e^{c_1c_2}_{j}e^{c_2}_{i}
=e^{c_2}_{j}e^{c_1c_2}_{i}e^{c_1}_{j}&
{\rm if }\,\,\lan \al^\vee_i,\al_j\ran
=\lan \al^\vee_j,\al_i\ran=-1,\\
&\hspace{-20pt}
e^{c_1}_{i}e^{c^2_1c_2}_{j}e^{c_1c_2}_{i}e^{c_2}_{j}
=e^{c_2}_{j}e^{c_1c_2}_{i}e^{c^2_1c_2}_{j}e^{c_1}_{i}&
{\rm if }\,\,\lan \al^\vee_i,\al_j\ran=-2,\,
\lan \al^\vee_j,\al_i\ran=-1,\\
&\hspace{-20pt}
e^{c_1}_{i}e^{c^2_1c_2}_{j}e^{c^3_1c_2}_{i}
e^{c^3_1c^2_2}_{j}e^{c_1c_2}_{i}e^{c_2}_{j}
=e^{c_2}_{j}e^{c_1c_2}_{i}e^{c^3_1c^2_2}_{j}e^{c^3_1c_2}_{i}
e^{c^2_1c_2}_je^{c_1}_i&
{\rm if }\,\,\lan \al^\vee_i,\al_j\ran=-3,\,
\lan \al^\vee_j,\al_i\ran=-1,
\end{array}
\]
\end{lem3}
{\sl Remark.} 
If $\lan\al^\vee_i,\al_j\ran\lan \al^\vee_j,\al_i\ran
\geq4$, there is no relation between $e_i$ and $e_j$.

\subsection{Unipotent Crystals}

In the sequel, we denote the unipotent subgroup 
$U^+$ by $U$. 
We define unipotent crystals (see \cite{BK}) associated to 
Kac-Moody groups. 


The definitions below are followed to \cite{BK}.
\begin{def3}
Let $X$ be an ind-variety over $\bbC$ and 
$\al:U\times X\rightarrow X$ be a rational $U$-action
such that $\al$ is defined on $\{e\}\times X$. Then, 
the pair ${\bf X}=(X,\al)$ is called a $U$-{\it variety}. 
For $U$-varieties ${\bf X}=(X,\al_X)$
and ${\bf Y}=(Y,\al_Y)$, 
a rational morphism
$f:X\rightarrow Y$ is called a 
$U$-{\it morphism} if it commutes with
the action of $U$.
\end{def3}
Now, we define the $U$-variety structure on $B^-=U^-T$.
By Proposition \ref{ind-g},
 $B^-$ is an ind-subgroup of $G$ and then
is an ind-variety over $\bbC$.
The multiplication map in $G$ induces the open embedding;
$ B^-\times U\hookrightarrow G,$
then this is a birational isomorphism. 
Let us denote the inverse birational isomorphism by $g$;
\[
 g:G\longrightarrow B^-\times U.
\]
Then we define the rational morphisms 
$\pi^-:G\rightarrow B^-$ and 
$\pi:G\rightarrow U$ by 
$\pi^-:={\rm proj}_{B^-}\circ g$ 
and $\pi:={\rm proj}_U\circ g$.
Now we define the rational $U$-action $\al_{B^-}$ on $B^-$ by 
\[
 \al_{B^-}:=\pi^-\circ m:U\times B^-\longrightarrow B^-,
\]
where $m$ is the multiplication map in $G$.
Then we obtain $U$-variety ${\bf B}^-=(B^-,\al_{B^-})$.
\begin{def3}
\label{uni-def}
\begin{enumerate}
\item
Let ${\bf X}=(X,\al)$ 
be a $U$-variety and $f:X \rightarrow {\bf B^-}$ 
be a $U$-morphism.
The pair $({\bf X}, f)$ is called 
a {\it unipotent $G$-crystal}
or, for short, {\it unipotent crystal}.
\item
Let $({\bf X},f_X)$ and $({\bf Y},f_Y)$ 
be unipotent crystals.
A $U$-morphism $g:X\ra Y$ is called a {\it morphism of 
unipotent crystals} if $f_X=f_Y\circ g$.
In particular, if $g$ is a birational isomorphism
of ind-varieties, it is called an {\it isomorphism of 
unipotent crystals}.
\end{enumerate}
\end{def3}
We define a product of 
unipotent crystlas following to \cite{BK}.
For unipotent crystals $({\bf X},f_X)$, $({\bf Y},f_Y)$, 
define a morphism 
$\al_{X\times Y}:U\tm X\tm Y\rightarrow X\tm Y$ by
\begin{equation}
\al_{X\tm Y}(u,x,y):=(\al_X(u,x),\al_Y(\pi(u\cdot f_X(x)),y)).
\label{XY}
\end{equation}
If there is no confusion,
we use abbreviated notation $u(x,y)$ 
for $\al_{X\tm Y}(u,x,y)$.
\begin{thm3}[\cite{BK}]
\label{prod}
\begin{enumerate}
\item
The morphism $\al_{X\tm Y}$ defined above 
is a rational $U$-morphism
on $X\tm Y$.
\item
Let ${\bf m}:B^-\tm B^-\rightarrow B^-$ 
be a multipication morphism 
and $f=f_{X\tm Y}:X\tm Y\rightarrow B^-$ be the 
rational morphism defined by 
\[
f_{X\tm Y}:={\bf m}\circ( f_X\tm f_Y).
\]
Then $ f_{X\tm Y}$ is a $U$-morphism and then, 
$({\bf X\tm Y}, f_{X\tm Y})$ is a unipotent crystal, 
which we call a product of unipotent crystals  
$({\bf X},f_X)$ and $({\bf Y},f_Y)$.
\item 
Product of unipotent crystals is associative.
\end{enumerate}
\end{thm3}
\subsection{From unipotent crystals to geometric crystals}
For $i\in I$, 
set $U^\pm_i:=U^\pm\cap \bar s_i U^\mp\bar s_i^{-1}$ and
$U_\pm^i:=U^\pm\cap \bar s_i U^\pm\bar s_i^{-1}$.
Indeed, $U^\pm_i=U_{\pm \al_i}$.
Set 
\[
 Y_{\pm\al_i}:=
\lan x_{\pm\al_i}(t)U_{\al}x_{\pm\al_i}(-t)
|t\in\bbC,\,\,\al\in \Delta^{\rm re}_{\pm}\setminus
\{\pm\al_i\}\ran.
\]
\begin{lem3}[\cite{Ku2},\cite{KP2}]
\label{dec-u}
For a simple root $\al_i$ $(i\in I)$, we have:
\begin{enumerate}
\item 
$Y_{\pm\al_i}=U^i_{\pm}$.
\item
$U^\pm=U_i^\pm\cdot Y_{\pm\al_i}$ (semi-direct product).
\item
$\bar s_i Y_{\pm\al_i}\bar s_i^{-1}=Y_{\pm\al_i}$.
\end{enumerate}
\end{lem3}
By this lemma, we have the unique 
decomposition;
\[
 U^-=U_i^-\cdot Y_{\pm\al_i}=U_{-\al_i}\cdot U^i_-.
\]
By using this decomposition, we get 
the canonical projection 
$\xi_i:U^-\rightarrow U_{-\al_i}$.
Now, we define the function on $U^-$ by 
\[
\chi_i:=y_i^{-1}\circ\xi_i:
U^-\longrightarrow U_{-\al_i}\longrightarrow
\bbC,
\]
and extend this to the function on 
$B^-$ by $\chi_i(u\cdot t):=\chi_i(u)$ for 
$u\in U^-$ and $t\in T$.
For a unipotent $G$-crystal $\bf(X,f_X)$, we define a function
$\vp_i:=\vp_i^X:X\rightarrow \bbC$ by 
\[
\vp_i:=\chi_i\circ{\bf f_X},
\]
and a ratinal morphism 
$\gamma_X:X\rightarrow T$ by 
\begin{equation}
\gamma_{X}:=
{\rm proj}_T\circ{\bf f_X}:X\rightarrow B^-\rightarrow T,
\label{gammax}
\end{equation}
where ${\rm proj}_T$ is the canonical projection.
Suppose that the function $\vp_i$ is not identially zero on $X$. 
We define a morphism $e_i:\bbC^\tm\tm X\rightarrow X$ by
\begin{equation}
e^c_i(x):=x_i
\left({\frac{c-1}{\vp_i(x)}}\right)(x).
\label{ei}
\end{equation}
\begin{thm3}[\cite{BK}]
\label{U-G}
For a unipotent $G$-crystal $\bf(X,f_X)$, 
suppose that 
the function $\vp_i$ is not identically zero
for any $i\in I$.
Then the rational morphisms $\gamma_X:X\rightarrow T$ 
and 
$e_i:\bbC^\tm\tm X\rightarrow X$ as above 
define a geometric 
$G$-crystal $(X,\gamma_X,\{e_i\}_{i\in I})$,
which is called the induced geometric $G$-crystals by 
unipotent $G$-crystal $({\bf X},f_X)$.
\end{thm3}
Note that in \cite{BK}, the cases $\vp_i\equiv0$ 
for some $i\in I$
are treated by considering Levi subgroups of $G$.
But here we do not treat such things.

The following product struture on geometric crystals are
most important results in the sense of comparison with
the tensor product theorem
in Kashiwara's crystal theory.
\begin{pro3}
For unipotent $G$-crystals $({\bf X},f_X)$
and $({\bf Y},f_Y)$, set 
the product $({\bf Z},f_Z):=({\bf X},f_X)
\tm({\bf Y},f_Y)$, where 
$Z=X\tm Y$. Let $(Z,\gamma_Z,\{e_i\})$ 
be the induced geometric 
$G$-crystal from $({\bf Z},f_Z)$.
Then we obtain;
\begin{enumerate}
\item $\gamma_Z={\bf m}\circ(\gamma_X\tm \gamma_Y)$.
\item For each $i\in I$, $(x,y)\in Z$, we obtain
\begin{equation}
\vp^Z_i(x,y)=
\vp^X_i(x)+\frac{\vp^Y_i(y)}{\al_i(\gamma_X(x))}.
\label{gamma-zxy}
\end{equation}
\item
For any $i\in I$, the action 
$e_i:\bbC^\tm\tm Z\rightarrow Z$ is given 
by: $e^c_i(x,y)=(e^{c_1}_i(x),e^{c_2}_i(y))$, where 
\begin{equation}
c_1=
\frac{c\al_i(\gamma_X(x))\vp^X_i(x)+\vp^Y_i(y)}
{\al_i(\gamma_X(x))\vp^X_i(x)+\vp^Y_i(y)},\,\,
c_2=
\frac{\al_i(\gamma_X(x))\vp^X_i(x)+\vp^Y_i(y)}
{\al_i(\gamma_X(x))\vp^X_i(x)+c^{-1}\vp^Y_i(y)}
\label{c1c2}
\end{equation}
\end{enumerate}
\end{pro3}
Here note that $c_1c_2=c$. 
The formula $c_1$ and $c_2$ in \cite{BK} 
seem to be different from ours.
Thus, we give the proof of (iii). 
Others are obtained by 
the same way as \cite{BK}.

{\sl Proof.}
By using the result (ii), we have
\[
\vp^Z_i(x,y)=
\vp^X_i(x)+\frac{\vp^Y_i(y)}{\al_i(\gamma_X(x))}.
\]
Here we set $A:=\frac{c-1}{\vp_i^Z(x,y)}$ for $(x,y)\in Z$.
Since by (\ref{XY}) we have 
\[
 \begin{array}{lll}
e^c_i(x,y)&=&x_i(A)(x,y)\\
&=& (x_i(A)(x),\pi(x_i(A)\cdot f_X(x))(y)),
\end{array}
\]
we get $\frac{c_1-1}{\vp^X_i(x)}=A$, and then we obtain 
$c_1$ in (\ref{c1c2}).

Let us see $c_2$. Writing $f_X(x)=u\cdot t$ 
($u\in U^-$, $t\in T$), 
by Lemma 3.1 (3.2) in \cite{BK}, we get
\[
\pi( x_i(A)\cdot f_X(x))
=x_i((A^{-1}+\chi_i(u)^{-1})^{-1}\al_i(t^{-1}))
\]
Since $\chi_i(u)=\vp_i(x)$ and 
$\al_i(t)=\al_i(\gamma_X(x))$, 
we obtain
\[
 \pi( x_i(A)\cdot f_X(x))=
x_i\left(
\frac{A}{(1+A\vp_i(x))\al_i(\gamma_X(x))}
\right).
\]
Now, set $B=\frac{A}{(1+A\vp^X_i(x))\al_i(\gamma_X(x))}$.
Substituting $A=\frac{c-1}{\vp_i^Z(x,y)}$  and 
$\frac{c_2-1}{\vp^Y_i(y)}=B$, we obtain the formula
$c_2$ in (\ref{c1c2}).\qed
\renewcommand{\thesection}{\arabic{section}}
\section{Crystal structure on Schubert varieties}
\setcounter{equation}{0}
\renewcommand{\theequation}{\thesection.\arabic{equation}}
\newtheorem{pro4}{Proposition}[section]
\newtheorem{thm4}[pro4]{Theorem}
\newtheorem{lem4}[pro4]{Lemma}
\newtheorem{ex4}[pro4]{Example}
\newtheorem{cor4}[pro4]{Corollary}
\theoremstyle{definition}
\newtheorem{def4}[pro4]{Definition}

\subsection{Highest weight modules and Schubert varieties}
As in Sect.2, let $G$ be a Kac-Moody group, 
$B^{\pm}=U^{\pm}T$ (resp. $U^{\pm}$)be the Borel 
(resp. unipotent)  subgroups in $G$ and 
$W$ be the associated Weyl group. Here, we have the
following Bruhat decomposition and Birkhoff decomposition;
\begin{pro4}[\cite{Ku2},\cite{PK},\cite{Slo}]
We have
\begin{eqnarray}
&&G=\bigcup_{w\in W}B^+\bar wB^+
=\bigcup_{w\in W}U^+\bar wB^+\q(\rm Bruhat \,\,decomposition),\\
&&G=\bigcup_{w\in W}B^-\bar wB^+
=\bigcup_{w\in W}U^-\bar wB^+\q(\rm Birkhoff\,\, decomposition).
\end{eqnarray}
\end{pro4}
Let $J\subset I$ be a subset of the index set $I$ and 
$W_J:=\lan s_i|i\in J\ran$ be the subgroup of $W$ associated 
with $J$. Set $P_J:=B^+ W_JB^+$ and call it a 
{\it $($standard$)$ parabolic subgroup} of 
$G$ associated with $J\subset I$. 
The following lemma is well-known;
\begin{lem4}
Any coset in $W/W_J$ contains a unique element $w^*$
of minimal length, and for any $w'\in W_J$, we have
$l(w^*w')=l(w^*)+l(w')$.
\label{min}
\end{lem4}
We denote the set of the elements $w^*$ 
as in Lemma \ref{min} by $W^J$, 
which is a set of representatives of $W/W_J$ in $W$.
There exist the following parabolic Bruhat/Birkhoff 
decompositions:
\begin{pro4}[\cite{Ku2},\cite{PK},\cite{Slo}]
\label{para-dec}
Let $J$ be a subset of $I$ and, $W_J$ and $W^J$ be as above. 
Then we have
\begin{eqnarray*}
&&G=\bigcup_{w^*\in W^J}U^+\bar{w^*}P_J, \\
&&G=\bigcup_{w^*\in W^J}U^-\bar{w^*}P_J.
\end{eqnarray*}
\end{pro4}

\subsection{Unipotent crystal structure on Schubert variety}
For $\Lm\in P_+$
($P_+$ is the set of dominant integral weight), let us denote
an integral highest weight simple module with
the highest weight $\Lm$ by $L(\Lm)$(\cite{Kac})
and its projective space by 
$\bbP(\Lm):=(L(\Lm)\setminus\{0\})/\bbC^\tm$.
Let $v_\Lm\in\bbP(\Lm)$ be the 
point correponding to the line
containing the highest weight vector of $L(\lm)$ and set 
\[
 X(\Lm):=G\cdot v_\Lm\subset \bbP(\Lm).
\]
Set $J_\Lm:=\{i\in I|\lan h_i,\Lm\ran=0\}$. By
Proposition \ref{para-dec} and the fact that $P_{J_\Lm}$ is the 
stabilizer of $v_\Lm$, we have the isomorphism between
$X(\Lm)$ and the flag variety $G/P_{J_\Lm}$:
\begin{pro4}[\cite{PK},\cite{Slo}]
\label{G/P}
There is the following isomorphism and the decomposition;
\[
\begin{array}{ccc}
\rho:
 G/P_{J_\Lm}
=\bigcup_{w\in W^{J_\Lm}}U^{\pm}\bar{w}P_{J_\Lm}/P_{J_\Lm}
&\mapright{\sim}& X(\Lm)\\
\qq\qq g\cdot P_{J_\Lm}& \mapsto& g\cdot v_\Lm
\end{array}
\]
\end{pro4}
\begin{def4}
We denote the image $\rho(U^+\bar{w}P_{J_\Lm}/P_{J_\Lm})$
(resp. $\rho(U^-\bar{w}P_{J_\Lm}/P_{J_\Lm})$)
by $X(\Lm)_w$ (resp. $X(\Lm)^w$) and call it
a {\it finite (resp. co-finite) Schubert 
cell} and its Zariski closure in 
$\bbP(\Lm)$ by $\ovl X(\Lm)_w$
(resp. $\ovl X(\Lm)^w$) and call
it a {\it finite (resp. co-finite) Schubert variety}.
\end{def4}
The names ``finite'' and ``co-finite'' come from the fact 
\[
 {\rm dim}X(\Lm)_w=l(w), \q {\rm codim}_{X(\Lm)}X(\Lm)^w=l(w),
\]
Indeed, $X(\Lm)_w\cong \bbC^{l(w)}$.
There exist the following closure relations;
\begin{equation}
\ovl X(\Lm)_w=\bigsqcup_{y\leq w, y\in W^{J_\Lm}}
X(\Lm)_y,\q
\ovl X(\Lm)^w=\bigsqcup_{y\geq w, y\in W^{J_\Lm}}X(\Lm)^y.
\end{equation}
Indeed, by \cite[7.1,7.3]{Ku2}, 
\begin{equation}
{\hbox{$\ovl X(\Lm)_w$ and 
$\ovl X(\Lm)^w$ are ind-varieties.}}
\label{ind}
\end{equation}
Let us associate a unipotent crystal structure
with $X(\Lm)_w$. Since by the definition of $X(\Lm)_w$ 
and Proposition \ref{G/P}, we have
$X(\Lm)_w=U^+\bar w\cdot v_\Lm$, which implies;
\begin{lem4}
\label{XU}
Schubert cell $X(\Lm)_w$ is a $U$-variety.
\end{lem4}
Next, let us construct a $U$-morphism 
$X(\Lm)_w\rightarrow B^-$.
For that purpose, we consider the following:
let 
$w=s_{i_1}s_{i_2}\cd s_{i_k}$ be a reduced expression and set
$U_w=U\cap \bar w U^-\bar w^{-1}$ and 
$U^w=U\cap \bar wU\bar w^{-1}$. Define
\[
\beta_1= \al_{i_1},\, \beta_2=s_{i_1}(\al_{i_2}),\,\cd,
\beta_k=s_{i_1}s_{i_2}\cd s_{i_{k-1}}(\al_{i_k}),
\]
then we have
\[
 U_w:=U_{\beta_1}\cdot U_{\beta_2}\cd U_{\beta_k}.
\]
This is a closed subgroup of $U$ and we have an
isomorphism of ind (algebraic)-varieties (\cite{Slo})
\begin{equation}
U_w\cong U_{\beta_1}\tm U_{\beta_2}\tm\cd \tm U_{\beta_k}
\cong \bbC^k,
\end{equation}
by 
\begin{eqnarray}
&& U_w\cdot \bar w =
U_{\al_{i_1}}\bar s_{i_1}\cdot
U_{\al_{i_2}}\bar s_{i_2}\cdot\cd
U_{\al_{i_k}}\bar s_{i_k}\mapright{\sim}
\bbC^k \label{cong}\\
&&x_{i_1}(a_1)\bar s_{i_1}\cdot
x_{i_2}(a_2)\bar s_{i_2}\cdot\cd
x_{i_k}(a_k)\bar s_{i_k}\mapsto
(a_1,a_2,\cd,a_k).\nn
\end{eqnarray}
\begin{lem4}[{\cite[2.2]{Slo}}]
\label{UU}
\begin{enumerate}
 \item We have a decomposition 
\begin{equation}
U=U_w\cdot U^w,
\label{U=UU}
\end{equation}
and the decomposition (\ref{U=UU}) 
is unique in the sense
;if $u_1v_1=u_2v_2$ $(u_i\in U_w,\,v_i\in U^w)$, then
$u_1=u_2$ and $v_1=v_2$.
\item
For any $w\in W^{J_\Lm}$ $(\Lm\in P_+)$, 
there exists an isomorphism of ind (algebraic)-varieties
\[
\begin{array}{ccc}
\delta:U_w &\mapright{\sim} &X(\Lm)_w \\
u&\mapsto &u\cdot v_\Lm
\end{array}
\label{UX}
\]
\end{enumerate}
\end{lem4}

The follwing lemma is the first step for our purpose.
\begin{lem4}
For any $u\in U$ and $w\in W$, there exist 
 unique
$u'\in U_w\cdot \bar w$ and $v\in U$ such that 
$u\bar w=u'v$.
\end{lem4}
{\sl Proof.}
By Lemma \ref{UU}(i), there are 
unique $u''\in U_w$ and $v''\in U^w$ such that
$u=u''v''$. By the definition 
$U^w=U\cap \bar wU\bar w^{-1}$, we have
$\bar w^{-1}v''\bar w\in U$. 
Thus, setting $u'=u''\bar w$ and 
$v=\bar w^{-1}v''\bar w$,
we get the desired result.\qed

By using this decomposition, we define the following 
rational morphisms;
\[
\begin{array}{cccc}
p_w:&U\cdot\bar w&\longrightarrow &U_w\cdot\bar w\\
&u \bar w&\mapsto& u'\\
p^w:&U\cdot\bar w&\longrightarrow &U\\
&u\bar w&\mapsto& v
\end{array}
\]
Define a rational $U$-action on $U_w\cdot\bar w$ by
\begin{eqnarray*}
U\tm U_w\cdot\bar w&\longrightarrow &U_w\cdot\bar w\\
(x,u\bar w)\,\,\,\,&\mapsto &x(u\bar w):=p_w(xu\bar w)
=xu\bar w\cdot p^w(xu\bar w)^{-1}
\end{eqnarray*}
Next, we show the following lemma:
\begin{lem4}
\label{alxuw}
Let $\pi^-:G\rightarrow B^-$ be as above. For $x\in U$ and 
$u\bar w\in U_w\bar w$, we have
$$
\al_{B^-}(x,\pi^-(u\bar w))=\pi^-(x(u\bar w)).
$$
\end{lem4}
\def\uw{u\bar w}
{\sl Proof.} We have
\[
 \begin{array}{ccl}
\pi^-(x(\uw))  &= &x(\uw)\cdot\pi(x(\uw))^{-1} \\
&=& x\uw\cdot p^w(x\uw)^{-1}
\pi(x\uw\cdot p^w(x\uw)^{-1})^{-1}\\
&=& x\uw\cdot p^w(x\uw)^{-1}p^w(x\uw)\pi(x\uw)^{-1} 
\q({\rm since} \,\,p^w(x\uw)\in U)\\
&=& x\uw\cdot \pi(x\uw)^{-1}=\pi^-(x\uw)
 \end{array}
\]
On the other hand, 
\[
 \begin{array}{ccl}
\al_{B^-}(x,\pi^-(\uw))
&= &\pi^-(x\pi^-(\uw))=x\pi^-(\uw)
\cdot\pi(x\pi^-(\uw))^{-1} \\
&=&x\uw\cdot\pi(\uw)^{-1}\cdot
\pi(x\uw\cdot\pi(\uw)^{-1})^{-1}\\
&=&x\uw\cdot\pi(\uw)^{-1}\cdot
\pi(\uw)\cdot\pi(x\uw)^{-1}
\q({\rm since} \,\,\pi(\uw)\in U)\\
&=&x\uw\cdot\pi(x\uw)^{-1}=\pi^-(x\uw),
 \end{array}
\]
which completes the proof.\qed

Define an isomorphism of ind (algebraic)-varieties 
\begin{eqnarray*}
\zeta:X(\Lm)_w&\mapright{\sim}& U_w\bar w\\
v&\mapsto & \zeta(v):=\delta^{-1}(v)\bar w,
\end{eqnarray*}
where $w\in W^{J_\Lm}$ and $\Lm\in P_+$.
Since $X(\Lm)_w$ is $U$-orbit of 
$\rho(\bar w\cdot P_{J_\Lm}/P_{J_\Lm})$, 
$U$ acts rationally on $X(\Lm)_w$. 
We denote the action of $x\in U$ on $v\in X(\Lm)_w$ by 
$x(v)$.
\begin{lem4}
The isomorphism $\zeta:X(\Lm)_w\rightarrow U_w\bar w$ is 
a $U$-morphism.
\end{lem4}
{\sl Proof.} It is sufficient to show that
$\zeta(x(v))=x(\zeta(v))$ for $x\in U$ and $v\in X(\Lm)_w$.
Set $u=\delta^{-1}(v)$ and then we have $v=u\bar wv_\Lm$.
Since $v_\Lm$ is stable by the action of $U$,{\it i.e.,}
$U\cdot v_\Lm=v_\Lm$, we get 
\[
 x(v)=p_w(x\uw)(v_\Lm).
\]
Since $p_w(x\uw)\in U_w\bar w$, we get
\[
 \zeta(x(v))=p_w(x\uw).
\]
We also have $x(\zeta(v))=x(\uw)=p_w(x\uw)$ and then 
$\zeta(x(v))=x(\zeta(v))$.\qed

\nd
Define a rational morphism $f_w:X(\Lm)_w\rightarrow B^-$ by
$f_w=\pi^-\circ\zeta$.
The following is one of the main results of this article.
\begin{thm4}
\label{uni-thm}
For $\Lm\in P_+$ and $w\in W^{J_\Lm}$, let $X(\Lm)_w$ 
be a finite Schubert cell and $f_w:X(\Lm)_w\rightarrow B^-$
be as defined above. Then the pair $(X(\Lm)_w,f_w)$ is a 
unipotent $G$-crystal.
\end{thm4}
{\sl Proof.}
We see that $X(\Lm)_w$ is a $U$-variety in Lemma \ref{XU}.
So, we may show that $f_w$ is a $U$-morphism.
For $x\in U$ and $v\in X(\Lm)_w$, we get 
\[
 f_w(x(v))=\pi^-(\zeta(x(v)))=
\pi^-(x(\zeta(v)))=\pi^-(x(\uw)),
\]
where $u=\delta^{-1}(v)$. On the other hand, 
\[
 x(f_w(v))=x(\pi^-(\zeta(v)))=x(\pi^-(\uw))
=\al_{B^-}(x,\pi^-(\uw)).
\]
By Lemma \ref{alxuw}, we obtain 
$ f_w(x(v))= x(f_w(v))$, which implies that
$f_w$ is a $U$-morphism.\qed

In the sense of Definition \ref{uni-def}(ii), $\zeta$
is an isomorphism of unipotent crystals on $X(\Lm)_w$
and $U_w \bar w$.

Since $X(\Lm)_w\hookrightarrow \ovl X(\Lm)_w $ 
is an open embedding, 
they are birationally equivalent. 
Let $\omega:\ovl X(\Lm)_w \rightarrow X(\Lm)_w$ be the
inverse birational isomorphism.
Thus, 
$\bar f_w:=f_w\circ\omega:\ovl X(\Lm)_w\rightarrow B^-$
is a $U$-morphism. Then we have
\begin{cor4}
Let $\ovl X(\Lm)_w$ be a finite Scubert variety and 
$\bar f_w$ be defined as above. Then the pair 
$(\ovl X(\Lm)_w, \bar f_w )$ is a unipotent $G$-crystal.
\end{cor4}

{\sl Remark.}
Note that for all $w\leq w'$, we have the closed 
embedding 
$\ovl X(\Lm)_w\hookrightarrow \ovl X(\Lm)_{w'}$
(\cite{Slo}), and then 
isomoprhism 
\[
X(\Lm) \,\,\,\mapright{\sim}
\lim_{\longrightarrow\atop{w\in W^{J_\Lm}}}
\ovl X(\Lm)_w.
\]
Nevertheless, in general, we do not obtain 
a unipotent crystal structure on $X(\Lm)$ by using 
this direct limit since for $y< w$, the rational 
morphism $\ovl f_w:\ovl X(\Lm)_w\ra B^-$ is not defined
on $\ovl X(\Lm)_y$.
\subsection{Geometric Crystal structure on $X(\Lm)_w$}

As we have seen in 3.3, we can associate 
geometric crystal structure with the finite Schubert cell
(resp. variety) $X(\Lm)_w$ (resp. $\ovl X(\Lm)_w$)
since we have shown that they are unipotent $G$-crystals.

Now, let us verify the condition by which the function 
$\vp_i:X(\Lm)_w\rightarrow \bbC$ is not identically zero.

We recall the formula:
\begin{equation}
x_i(a)y_j(b)=
\left\{
\begin{array}{ll}
 y_i(\frac{b}{1+ab})\al_i^\vee(1+ab)x_i(\frac{a}{1+ab}) &
{\rm if }\,\,i=j\\
 y_j(b)x_i(a)&{\rm if}\,\,i\ne j
\end{array}
\right.
\label{xy=yax}
\end{equation}
Hence, we have
\begin{equation}
x_i(c)\bar s_i=y_i(\frac{1}{c})
\al_i^\vee(c)x_i(-\frac{1}{c}),
\label{xs}
\end{equation}
where $\bar s_i=x_i(-1)y_i(1)x_i(-1)$.
We also have
\begin{equation}
\al_i^\vee(a)x_j(b)=x_j(a^{a_{ij}}b)
\al_i^\vee(a),
\q
\al_i^\vee(a)y_j(b)=y_j(a^{-a_{ij}}b)
\al_i^\vee(a)
\label{ae}
\end{equation}
By the formula (\ref{xy=yax}), (\ref{xs}) 
and (\ref{ae}), we obtain 
\begin{eqnarray}
x_i(a)\cdot\left(y_j(\frac{1}{c})\al_j^\vee(c)\right)
&=& \left(y_j(\frac{1}{c})\al_j^\vee(c)\right)
\cdot x_i(c^{-a_{ji}}a),\q(i\ne j)
\label{xyij}\\
x_i(a)\cdot\left(y_i(\frac{1}{c})\al_i^\vee(c)\right)
&=& \left(y_i(\frac{1}{a+c})\al_i^\vee(a+c)\right)
\cdot x_i(\frac{a}{ac+c^2}).
\label{xyii}
\end{eqnarray}
Due to these formula, we get the following lemma;
\begin{lem4}
\label{xxyy}
For $w=s_{i_1}s_{i_2}\cd s_{i_k}\in W$ 
(reduced expression) and 
$c_1,c_2,\cd,c_k\in \bbC^\tm$, 
there exist $c'_1,c'_2,\cd,c'_k$
such that
\begin{eqnarray}
&&\pi^-(x_{i_1}(c_1)\bar s_{i_1}\cdot
x_{i_2}(c_2)\bar s_{i_2}\cd
x_{i_k}(c_k)\bar s_{i_k})\nn\\
&&\qq\qq=y_{i_1}(\frac{1}{c'_1})\al_{i_1}^\vee(c'_1)\cdot
y_{i_2}(\frac{1}{c'_2})\al_{i_2}^\vee(c'_2)\cd
y_{i_k}(\frac{1}{c'_k})\al_{i_k}^\vee(c'_k)
\label{im-pi}
\end{eqnarray}
\end{lem4}
For $w\in W$, let $w=s_{i_1}s_{i_2}\cd s_{i_k}$ be 
one reduced expression and set 
\[
 I(w):=\{i_1,i_2,\cd,i_k\}.
\]
Indeed, this does ont depend on 
the choice of reduced expressions
since $W$ is a Coxeter group.
By Lemma \ref{xxyy}, we get
\begin{lem4}
For $w\in W$ and $i\in I$, if $i\in I(w)$, then 
the function 
$\vp_i:X(\Lm)_w\rightarrow \bbC$ is not identically 
zero.
\end{lem4}
Now, by Theorem \ref{U-G}, we have
\begin{thm4}
For $w\in W$, suppose that $I=I(w)$. We can associate the 
geometric $G$-crystal structure with the finite Schubert 
cell $X(\Lm)_w$ $($resp. variety $\ovl X(\Lm)_w$  $)$ 
by setting
$($see $(\ref{gammax})$ and $(\ref{ei})$$)$
\[
 \gamma_w:={\rm proj}_T\circ f_w\,\,({\rm resp.} \,
\ovl\gamma_w:={\rm proj}_T\circ \bar f_w),\q 
e^c_i(x)=x_i\left(\frac{c-1}{\vp_i(x)}\right)(x),
\]
where ${\rm proj}_T:B^-=U^-T\ra T$.
\end{thm4}
We denote this induced geometric crystal by 
$(X(\Lm)_w,\gamma_w,\{e_i\}_{i\in I})$ 
(resp. $(\ovl X(\Lm)_w,\ovl\gamma_w,
\{e_i\}_{i\in I})$).
This geometric/unipotent crystal 
$(X(\Lm)_w,\gamma_w,\{e_i\}_{i\in I})$ is realized in 
$B^-$ in the following sense.
\begin{pro4}
For $w=s_{i_1}\cd s_{i_k}$, define 
\[
 B^-_w:=\{Y_w(c_1,\cd,c_k):=
y_{i_1}(\frac{1}{c_1})\al_{i_1}^\vee(c_{i_1})
\cd y_{i_k}(\frac{1}{c_k})\al_{i_k}^\vee(c_{i_k})
\in B^-
|c_i\in \bbC^\tm\}.
\]
and $U$-actions on $B^-_w$ by 
\[
 u(Y_w(c_1,\cd,c_k)):=\pi^-(u\cdot Y_w(c_1,\cd,c_k))
\q(u\in U).
\]						     
Then $X(\Lm)_w$ and $B^-_w$ are 
birationally equivalent
via $f_w$ and isomorphic as unipotent crystals.
Moreover, they are isomorphic 
as induced geometric crystals.
\end{pro4}
{\sl Proof.} By Lemma \ref{xxyy},it is sufficient to show 
that they are birationally equivalent to each other
and then we may show that
$U_w\cdot w$ and $B^-_w$ are birationally equivalent
via $\pi^-$.
For that purpose,  since we have 
the isomorphism (\ref{cong}) and 
the birational isomorphism $B^-_w\cong (\bbC^\tm)^k$,
 it sufficies to show that 
the correpondence $(c_1,\cd,c_k)\longleftrightarrow
(c'_1,\cd,c'_k)$ in (\ref{im-pi}) is birational.
In (\ref{im-pi}),
each $c'_i$ is a rational function 
in $c_1,c_2,\cd,c_i$ 
obtained by composing the birational morphisms defined
by (\ref{xyij}) and (\ref{xyii})
(in particular, $c_1=c'_1$),
which implies that
$U_w\cdot w$ and $B^-_w$ 
are birationally equivalent.\qed

\begin{ex4}
\label{ex-sym}
We consider the case $G=SL_{n+1}(\bbC)$,{\it i.e.,}
the Cartan matrix $A=(a_{ij})_{i,j\in I}$ is given by;
$a_{ii}=2$, $a_{ii+1}=-1$ and $a_{ij}=0$ otherwise.
Here $I=\{1,2,\cd,n\}$.
Take $\til w=s_1s_2\cd s_n\in W$. In this cse, 
we can easily find 
that $I=I(\til w)$ and 
\[
 \pi^-(x_1(c_1)\bar s_1x_2(c_2)\bar s_2
\cd x_n(c_n)\bar s_n)
=y_1(\frac{1}{c_1})\al^\vee_1(c_1)
y_2(\frac{1}{c_2})\al^\vee_2(c_2)
\cd y_n(\frac{1}{c_n})\al^\vee_n(c_n).
\]
Here changing the coordinate by $c_i=a_1a_2\cd a_i$ and 
identifying 
$y_i(a)=I_n+aE_{i+1\,i}$, we obtain
\[
 f_{\til w}(X(\Lm)_{\til w})=
\left\{u(a):=\left(
\begin{array}{ccccccc}
a_1 & & & & & &\\
1   &a_2 & & & & &\\
 &1 & \cdot& & & &\\
 & &\cdot &\cdot & & &\\
 & & &\cdot &\cdot & &\\
 & & & & &a_n & \\
 & & & & &1 &\frac{1}{a_1\cd a_n} 
\end{array}
\right);a_i\in \bbC^\tm
\right\}
\]
where $a=(a_1,\cd,a_{n+1})$ and $a_1a_2\cd a_{n+1}=1$.
By using this explicit presentation, we describe the
geometric crystal structure of 
$f_{\til w}(X(\Lm)_{\til w})$.
Since $\vp_i(u(a))=\frac{1}{a_i}$, we have
\[
 e_i^c(u(a))=x_i(a_i(c-1))\cdot u(a)
\cdot x_i(a_{i+1}(c^{-1}-1))
=u(a_1,\cd,ca_i,c^{-1}a_{i+1},\cd,a_{n+1}).
\]
Further, we have $\gamma_{\til w}$ by
\[
 \gamma_{\til w}
(x_1(c_1)\bar s_1x_2(c_2)\bar s_2\cd x_n(c_n)\bar s_n)
=\al^\vee_1(c_1)\al^\vee_2(c_2)\cd \al^\vee_n(c_n).
\]
\end{ex4}

\renewcommand{\thesection}{\arabic{section}}
\section{Tropicalization of 
Crystals and Schubert Varieties}
\setcounter{equation}{0}
\renewcommand{\theequation}{\thesection.\arabic{equation}}
\theoremstyle{definition}
\newtheorem{def5}{Definition}[section]
\theoremstyle{plain}
\newtheorem{thm5}[def5]{Theorem}
\newtheorem{lem5}[def5]{Lemma}
\newtheorem{ex5}[def5]{Example}
\newtheorem{pro5}[def5]{Proposition}
\newtheorem{cor5}[def5]{Corollary}

We use the same notations 
as in the previous sections otherwise stated.
We introduce a positive structure 
on geometric crystals and 
their ultra-discritizations and tropicalizations following
to \cite[2.5]{BK}.

Let $T$ be an algebraic torus over $\bbC$ and 
$X^*(T)$ (resp. $X_*(T)$) be the lattice of characters
(resp. co-characters)
of $T$. Let $R:=\bbC[[c]][c^{-1}]$ and set
$L(T):=\{\phi\in {\rm Hom}(O_T,R)\}$ 
($O_T$ is the ring of 
regular functions on $T$), which is called 
a set of {\it formal loops} on $T$.
Here we specify the  discrete valuation 
$$
\begin{array}{cccc}
v:&R\setminus\{0\}&\longrightarrow &\ZZ\\
&\sum_{n>-\ify} a_n c^n&\mapsto
&-{\rm min}\{n\in\ZZ|a_n\ne0\}.
\end{array}
$$
For any $\phi\in L(T)$, 
set 
${\rm deg}_T(\phi):=v\circ\phi|_{X^*(T)}$. Since 
for $f_1,f_2\in R\setminus\{0\}$
\begin{equation}
v(f_1 f_2)=v(f_1)+v(f_2),
\label{ff=f+f}
\end{equation}
${\rm deg}_T(\phi)$ can be 
considered as an element in 
$X_*(T)={\rm Hom}(X^*(T),\ZZ)$.
Hence, ${\rm deg}_T$ can be seen as a map 
${\rm deg}_T:L(T)\ra X_*(T)$.
For any $\lm^\vee\in X_*(T)$, define
$L_{\lm^\vee}(T)
:={\rm deg}_T^{-1}(\lm^\vee)\subset L(T)$.
Since ${\rm deg}_T^{-1}(\lm^\vee)$ 
has an irreducible 
pro-$\bbC$ variety structure
and $L(T)=\bigsqcup_{\lm^\vee\in X_*(T)}
L_{\lm^\vee}(T)$, 
the set of irreducible components 
$\pi_0(L(T))=\{L_{\lm^\vee}(T)|\lm^\vee\in X_*(T)\}$ 
can be identified with $X_*(T)$, {\it i.e.,}
${\rm deg}_T$ induces the bijection 
$\til{\rm deg}_T:\pi_0(L(T))\mapright{\rm 1:1}X_*(T)$.

More explicitly, set $T=(\bbC^\tm)^l$ and  identify 
$L(T)$ with $(R^\tm)^l$.
For $\lm^\vee(c)=(c^{m_1},c^{m_2},\cd,c^{m_l})$
($m_j\in \ZZ$), we have
$$
L_{\lm^\vee}(T)
=\left\{\left(b_1 c^{-m_1}+\sum_{n>-m_1}a_nc^n,\cd,
b_l c^{-m_l}+\sum_{n>-m_l}a_nc^n\right):
b_1,\cd,b_l\ne0\right\}.
$$

Let $f:T\rightarrow T'$ be a rational morphism 
between two algebaric tori $T$ and $T'$.
The morphism $f$ induces the rational morphism 
$\til f:L(T)\rightarrow L(T')$ and then the map
$\pi_0(\til f):\pi_0(L(T))\rightarrow \pi_0(L(T'))$,
which defines the map 
${\rm deg}(f):X_*(T)\rightarrow X_*(T')$.
\[
\xymatrix@R=3pt{
\pi_0(L(T))
\ar@{->}^{\pi_0(\til f)} @<1pt>[rr]
\ar@<-1pt>[ddd]^{\til{\rm deg}_T}
&& \pi_0(L(T'))
\ar@<-1pt>[ddd]^{\til{\rm deg}_{T'}}\\ \\ \\
X_*(T) \ar@{->}[rr]^{{\rm deg}(f)}&& X_*(T')
}
\]

A rational function $f(c)\in \bbC(c)$ $(f\ne0)$ is 
{\it positive} if $f$ can be expressed as a ratio
of polynomials with positive coefficients.

\nd
{\sl Remark.}
A rational function $f(c)\in \bbC(c)$ is positive 
if and only if $f(a)>0$ for any $a>0$
(pointed out by M.Kashiwara).

If $f_1,\,\,f_2\in \bbC(c)(\subset R)$ 
are positive, then we have
\begin{eqnarray}
&& v(f_1f_2)=v(f_1)+v(f_2), \\
&& v\left(\frac{f_1}{f_2}\right)=v(f_1)-v(f_2),\\
&& v(f_1+f_2)={\rm max}(v(f_1),v(f_2)).
\end{eqnarray}

\begin{def5}[\cite{BK}]
A rational morphism $f:T\rightarrow T'$ between
two algebraic tori $T,T'$ is called {\it positive},
if the following two conditions are satisfied:
\begin{enumerate}
\item For any co-character 
$\lm^\vee:\bbC^\tm\rightarrow T$, the image of $\lm^\vee$
is contained in dom$(f)$.
\item For any co-character 
$\lm^\vee:\bbC^\tm\rightarrow T$ and 
any character 
$\mu:T'\rightarrow \bbC^\tm$, the composition 
$\mu\circ f\circ \lm^\vee$ is a positive rational 
function.
\end{enumerate}
\end{def5}

Denote by ${\rm Mor}^+(T,T')$ the set of 
positive rational morphisms from $T$ to $T'$.

\begin{lem5}[\cite{BK}]
For any positive rational morphisms 
$f\in {\rm Mor}^+(T_1,T_2)$             
and $g\in {\rm Mor}^+(T_2,T_3)$, 
 the composition $g\circ f$
is in ${\rm Mor}^+(T_1,T_3)$.
\end{lem5}
By this lemma, we can define a category ${\cal T}_+$
whose objects are algebraic tori over $\bbC$ and arrows
are positive rational morphisms.
\begin{lem5}[\cite{BK}]
For any algebraic tori $T_1$, $T_2$, $T_3$, 
and positive rational morphisms 
$f\in {\rm Mor}^+(T_1,T_2)$, 
$g\in {\rm Mor}^+(T_2,T_3)$, we have
\[
 {\rm deg}(g\circ f)={\rm deg}(g)\circ{\rm deg}(f).
\]
\end{lem5}
By this lemma, we obtain a functor 
\[
\begin{array}{cccc}
{\cal UD}:&{\cal T}_+&\longrightarrow &{{\hbox{\germ Set}}}\\
&T&\mapsto& X_*(T)\\
&(f:T\rightarrow T')&\mapsto& 
({\rm deg}(f):X_*(T)\rightarrow X_*(T')))
\end{array}
\]


\begin{def5}[\cite{BK}]
\begin{enumerate}
\item
Let $\chi=(X,\gamma,\{e_i\}_{i\in I})$ be a 
geometric pre-crystal, $T$ be an algebraic torus
and $\theta:T\rightarrow X$ 
be a birational isomorphism.
The isomorphism $\theta$ is called 
{\it positive structure} on
$\chi$ if it satisfies
\begin{enumerate}
\item the rational morphism 
$\gamma\circ \theta:T'\rightarrow T$ is positive.
\item
For any $i\in I$, the rational morphism 
$e_{i,\theta}:\bbC^\tm \tm T'\rightarrow T'$ given by
\[
e_{i,\theta}(c,t)
:=\theta^{-1}\circ e_i^c\circ \theta(t)
\]
is positive.
\end{enumerate}
\item
Two positive structures $\theta_1,\theta_2$ 
on a geometric pre-crystal are {\it equivalent} if 
the rational morphisms $\theta_1^{-1}\circ 
\theta_2$ and 
$\theta_2^{-1}\circ \theta_1$ are positive.
\end{enumerate}
\end{def5}
Applying the functor ${\cal UD}$ 
to positive rational morphisms
$e_{i,\theta}:\bbC^\tm \tm T'\rightarrow T'$ and
$\gamma\circ \theta:T'\ra T$
(the notations are
as above), we obtain
\begin{eqnarray*}
\til e_i&:=&{\cal UD}(e_{i,\theta}):
\ZZ\tm X_*(T) \rightarrow X_*(T)\\
\til\gamma&:=&{\cal UD}(\gamma\circ\theta):
X_*(T')\rightarrow X_*(T).
\end{eqnarray*}
Now, for given positive structure $\theta:T'\rightarrow X$
on a geometric pre-crystal 
$\chi=(X,\gamma,\{e_i\}_{i\in I})$, we associate 
the triplet $(X_*(T'),\til \gamma,\{\til e_i\}_{i\in I})$
with a pre-crystal structure (see \cite[2.2]{BK}) 
and denote it by ${\cal UD}_{\theta,T'}(\chi)$.
By Lemma \ref{Verma}, we have the following theorem:

\begin{thm5}
For any geometric crystal 
$\chi=(X,\gamma,\{e_i\}_{i\in I})$ and positive structure
$\theta:T'\rightarrow X$, the associated pre-crystal 
${\cal UD}_{\theta,T'}(\chi)=
(X_*(T'),\til\gamma,\{\til e_i\}_{i\in I})$ 
is a free $W$-crystal (see \cite[2.2]{BK})
\end{thm5}

We call the functor $\cal UD$
{\it ``ultra-discritization''}
instead of ``tropicalization'' unlike in \cite{BK}.
The term ``tropicalization'' here means the inverse
tropicalization in \cite{BK}.
More precisely, 
for an object $B$ in ${{\hbox{\germ Set}}}$, if there
exists an object $T$ in ${\cal T}_+$ such that 
${\cal UD}(T)\cong B$ as crystals, 
we call $T$
a {\it tropicalization} of $B$.

Now, we define certain positive strucuture on 
geometric crystal $B^-_w$ 
$(I=I(w),{\hbox{ and }}w\in W^{J_\Lm})$
and see that it turns out to be 
a tropicalization of (Langlands dual of) 
some Kashiwara's crystal.

For $i\in I$, let $B_i$ be the crystal defined by 
(see {\it e.g.}\cite{K3})
\begin{eqnarray*}
& B_i:=\{(x)_i|x\in \ZZ\},\\
&\til e_i(x)_i=(x+1)_i, \,\,
\til f_i(x)_i=(x-1)_i, \,\, 
\til e_j(x)_i=\til f_j(x)_i=0\,\,(i\ne j)\\
&wt(x)_i=x\al_i, \vep_i(x)_i=-x,\,\,\vp_i(x)_i=x,\,\,
\vep_j(x)_i=\vp_j(x)_i=-\ify\,\,(i\ne j).
\end{eqnarray*}

For $w=s_{i_1}s_{i_2}\cd s_{i_k}\in W$
and ${\bf i}=(i_1,i_2,\cd,i_k)\in R(w)$, 
we define the morphism
$\theta_{\bf i}:(\bbC^\tm)^k\rightarrow B^-_w$ by 
\begin{equation}
 \theta_{\bf i}(c_1,c_2,\cd,c_k)
:=y_{i_1}(\frac{1}{c_1})\al_{i_1}^\vee(c_1)
\cd y_{i_k}(\frac{1}{c_k})\al_{i_k}^\vee(c_k)
\label{posi}
\end{equation}
Similar statements to the following proposition 
are given in \cite[Theorem 2.11]{BK}
for reductive cases. Here we show it
for arbitrary Kac-Moody cases
by direct methods.
\begin{pro5}
\label{BandB}
\begin{enumerate}
 \item 
For any ${\bf i}\in R(w)$ $(w\in W)$, 
the morphism $\theta_{\bf i}$ defiend in (\ref{posi})
is a positive structure on the geometric crystal
$B^-_w$.
\item
Geometric crystal $B^-_w$ is a tropicalization of
the Langlands dual of 
the crystal $B_{i_1}\ot B_{i_2}\ot\cd\ot B_{i_k}$
with respect to  the positive structure 
$\theta_{\bf i}(c_1,c_2,\cd,c_k)$, or
equivalently 
${\cal UD}(B^-_w)\cong 
{\rm Langlands\,\,dual}(B_{i_1}\ot\cd\ot B_{i_k})$ as crystals.
\end{enumerate}
\end{pro5}

{\sl Proof.}
It is clear that $\theta_{\bf i}$ 
is a birational isomorphism. 
Since the ratinal morphism 
$\gamma:B^-_w\rightarrow T$ is given 
by
\[
 \gamma\left(
y_{i_1}(\frac{1}{c_1})\al_{i_1}^\vee(c_1)
\cd y_{i_k}(\frac{1}{c_k})\al_{i_k}^\vee(c_k)\right)
=\al_{i_1}^\vee(c_1)\cd \al_{i_k}^\vee(c_k),
\]
we have that $\gamma\circ\theta_{\bf i}$ is positive.
In order to show that 
$e_{i,\theta_{\bf i}}:\bbC^\tm \tm T' \ra T'$ is positive,
we see the explicit action of 
$e_i^c$ on $Y_w(c_1,\cd,c_k)$.
First let us evaluate $\vp_i(Y_w(c_1,\cd,c_k))$. 
\begin{lem5}
\label{vp-Y}
 For $Y:=y_{i_1}(a_1)\cd y_{i_k}(a_k)\in U^-$, we have
\begin{equation}
\vp_i(Y)
=\sum_{i_j=i}a_{i_j}.
\label{vp-y}
\end{equation}
\end{lem5}
{\sl Proof.}
Let $\{j_1,j_2,\cd,j_r\}$ 
$(j_1<j_2<\cd <j_r)$be the set of indecies such that
$i_{j_m}=i$. Then we can write
\[
Y=A_0 \cdot y_i(a_{i_{j_1}}) \cdot A_1
\cdot y_i(a_{i_{j_2}})\cdot A_2\cdot
y_i(a_{i_{j_3}})
\cd A_{r-1}\cdot y_i(a_{i_{j_r}})\cdot A_r,
\]
where $A_s:=\prod_{j_s<p<j_{s+1}}y_{i_p}(a_{i_p})$ 
$(j_0=0,j_{r+1}=k+1)$.
Here we set
\[
 B_m:=y_i(-\sum_{m<s\leq r}a_{i_{j_s}})\cdot 
A_m \cdot y_i(\sum_{m<s\leq r}a_{i_{j_s}}),
\]
Then we have
\begin{equation}
 Y=y_i(\sum_{0<s\leq r}a_{i_{j_s}})\cdot 
(B_0\cdot B_1\cd B_r).
\label{Y-dec}
\end{equation}
Since $B_0\cdot B_1\cd B_r$ is in $Y_{-\al_i}$ and 
the decomposition (\ref{Y-dec})
is unique by Lemma \ref{dec-u}, we have 
\[
 \vp_i(Y)=\sum_{0<s\leq r}a_{i_{j_s}}
=\sum_{i_j=i}a_{i_j},
\]
which is the desired result.\qed

Set
\[
C^*_j:=(c_1^{a_{i_1,i_j}}c_2^{a_{i_2,i_j}}\cd 
c_{j-1}^{a_{i_{j-1},i_j}}c_j)^{-1}\,\,\qq(C^*_1=1/c_1),
\]
where $a_{i,j}$ is an $(i,j)$-entry of 
the generalized Cartan matrix $A$.
By (\ref{ae}), we have 
\[
Y_w(c_1,\cd, c_k)
=y_{i_1}(C^*_1)\cd y_{i_k}(C^*_k)\al^\vee_{i_1}(c_1)
\cd \al^\vee_{i_k}(c_k).
\]
Then by Lemma \ref{vp-Y}, we obtain
\begin{equation}
\vp_i(Y_w(c_1,\cd, c_k))
=\sum_{i_j=i} C^*_{i_j}
=\sum_{j=1}^k
\frac{\delta_{i,i_j}}{c_1^{a_{i_1,i}}
c_2^{a_{i_2,i}}\cd 
c_{j-1}^{a_{i_{j-1},i}}c_j}
\label{vp-Yw}
\end{equation}

For $c\in \bbC$ and $i\in I$, 
define $\{\ovl C_j\}_{1\leq j\leq k}$ and 
$\{\til C_j\}_{0\leq j\leq k}$ recursively by 
\[
\ovl C_0=c, \qq
\til C_j={c_j+\del_{i_j,i}\ovl C_{j-1}},\qq
 \ovl C_j=\frac{\ovl C_{j-1}\cdot c_j^{1-a_{i_j,i}}}
{\til C_j}.
\]
Then, by using (\ref{xyij}) and (\ref{xyii})
repeatedly,
we obtain 
\begin{equation}
x_i(c)(Y_w(c_1,\cd,c_k))
=Y_w(\til C_1,\cd,\til C_k).
\label{xi-Y}
\end{equation}

It is easy to get the explicit form of $\ovl C_j$:
\[
 \ovl C_j=
\frac{\displaystyle c\prod_{m=1}^j c_m^{1-a_{i_m,i}}}
{\displaystyle\sum_{1\leq m\leq j,i_m=i} c\cdot D_m
+\prod_{m=1}^j c_m},
\]
where 
\[
D_m:={c_1^{1-a_{i_1,i}}\cd c_{m-1}^{1-a_{i_{m-1},i}}
\cdot c_{m+1}\cd c_{j-1}c_j}.
\]
Now, in (\ref{xi-Y}) 
replacing $c$ with $(c-1)/\vp_i(Y_w(c_1,\cd,c_k))$ and 
using (\ref{vp-Yw}),
we obtain 
\[
e_i^c(Y_w(c_1,\cd,c_k))
=x_i\left(\frac{c-1}{\vp_i(Y_w(c_1,\cd,c_k))}\right)
(Y_w(c_1,\cd,c_k)))
=:Y_w({\cal C}_1,\cd,{\cal C}_k),
\]
where
\begin{equation}
{\cal C}_j:=
c_j\cdot \frac{\displaystyle \sum_{1\leq m\leq j,i_m=i}
 \frac{c}
{c_1^{a_{i_1,i}}\cd c_{m-1}^{a_{i_{m-1},i}}c_m}
+\sum_{j< m\leq k,i_m=i} \frac{1}
{c_1^{a_{i_1,i}}\cd c_{m-1}^{a_{i_{m-1},i}}c_m}}
{\displaystyle\sum_{1\leq m<j,i_m=i} 
 \frac{c}
{c_1^{a_{i_1,i}}\cd c_{m-1}^{a_{i_{m-1},i}}c_m}+
\mathop\sum_{j\leq m\leq k,i_m=i}  \frac{1}
{c_1^{a_{i_1,i}}\cd c_{m-1}^{a_{i_{m-1},i}}c_m}}.
\label{eici}
\end{equation}
By this formula, it is clear that 
$e_{i,\theta_{\bf i}}$ is positive. We have shown (i).

Next, in order to show (ii), 
we see the action of $\til e_i^c$
on $B_{i_1}\ot \cd\ot B_{i_k}$.
Take $b_{\bf i}=(b_1)_{i_1}\ot\cd\ot (b_k)_{i_k}$
$({\bf i}=(i_1,\cd,i_k),\,\,b_j\in \ZZ)$.
Since the action of $\til e_i$ on 
tensor products is described
explicitly in \cite{K3}, we obtain 
\[
\til e_i^c(b_{\bf i})
=(\beta_1)_{i_1}\ot\cd\ot (\beta_k)_{i_k},
\]
where
\begin{eqnarray}
\beta_j= b_j+{\rm max}\left(
\mathop{\rm max}_{\begin{array}{c}\scriptstyle
1\leq m\leq j, \\ \scriptstyle i_m=i\end{array}}
(c-b_m-\sum_{l<m}b_l a_{i,i_l}),\,\,
\mathop{\rm max}_{\begin{array}{c}\scriptstyle
j< m\leq k,\\ \scriptstyle i_m=i\end{array}}
(-b_m-\sum_{l<m}b_l a_{i,i_l})\right) \nn  \\
-{\rm max}
\left(\mathop{\rm max}_{\begin{array}{c}\scriptstyle
1\leq m< j,\\ \scriptstyle i_m=i\end{array}}
(c-b_m-\sum_{l<m}b_l a_{i,i_l}),\,\,
\mathop{\rm max}_{\begin{array}{c}\scriptstyle
j\leq m\leq k,\\ \scriptstyle i_m=i\end{array}}
(-b_m-\sum_{l<m}b_l a_{i,i_l})\right)\qq\qq
\label{eibi}
\end{eqnarray}
Now, we know that (\ref{eici}) and (\ref{eibi})
are related to each other by the 
tropicalization/ultra-discritzation operations:
\[
\xymatrix@R=3pt{
{\cal C}_j\ar@{->}^{
\rm\scriptstyle ultra-discritization} @<1pt>[rrr]
\ar@{<-}_{
\rm\scriptstyle tropicalization} @<-1pt>[rrr]
&&& \beta_j \\
c_j \ar@{<->}[rrr]&&& b_j\\
x\cdot y \ar@{<->}[rrr]&&& x+y \\
\frac{x}{y} \ar@{<->}[rrr]&&& x-y \\
x+y \ar@{<->}[rrr]&&& {\rm max}(x,y)\\
a_{i,j} \ar@{<->}_{\rm \scriptstyle
Langlands\,\, dual}[rrr]&&& a_{j,i}
}
\]
We have completed the proof of (ii).\qed

The formula similar to (\ref{eici}), (\ref{eibi})
are given in \cite[5.2.]{BZ} for the longest element
$w_{0}$ (in reductive cases).

The following formulae are an immediate consequence of 
Proposition \ref{BandB} and  Lemma \ref{Verma}, which 
are given implicitly in \cite{K3}
and shown by direct method in 
\cite{NS}.
\begin{cor5}
On the crystal $B_{i_1}\ot \cd\ot B_{i_k}$, we have
for any $c_1,c_2\in \ZZ_{\geq0}$
\[
 \begin{array}{lll}
&\hspace{-20pt}\til e^{c_1}_{i}\til e^{c_2}_{j}
=\til e^{c_2}_{j}\til e^{c_1}_{i}&
{\rm if }\,\,\lan \al^\vee_i,\al_j\ran=0,\\
&\hspace{-20pt} \til e^{c_1}_{i}
\til e^{c_1+c_2}_{j}\til e^{c_2}_{i}
=\til e^{c_2}_{j}\til e^{c_1+c_2}_{i}\til e^{c_1}_{j}&
{\rm if }\,\,\lan \al^\vee_i,\al_j\ran
=\lan \al^\vee_j,\al_i\ran=-1,\\
&\hspace{-20pt}
\til e^{c_1}_{i}\til e^{2c_1+c_2}_{j}
\til e^{c_1+c_2}_{i}\til e^{c_2}_{j}
=\til e^{c_2}_{j}\til e^{c_1+c_2}_{i}
\til e^{2c_1+c_2}_{j}\til e^{c_1}_{i}&
{\rm if }\,\,\lan \al^\vee_i,\al_j\ran=-1,\,
\lan \al^\vee_j,\al_i\ran=-2,\\
&\hspace{-20pt}
\til e^{c_1}_{i}\til e^{2c_1+c_2}_{j}\til e^{3c_1+c_2}_{i}
\til e^{3c_1+2c_2}_{j}
\til e^{c_1+c_2}_{i}\til e^{c_2}_{j}
&\\
&\qq\q =\til e^{c_2}_{j}\til e^{c_1+c_2}_{i}
\til e^{3c_1+2c_2}_{j}
\til e^{3c_1+c_2}_{i}
\til e^{2c_1+c_2}_je^{c_1}_i&
{\rm if }\,\,\lan \al^\vee_i,\al_j\ran=-1,\,
\lan \al^\vee_j,\al_i\ran=-3.
\end{array}
\]
\end{cor5}
{\sl Remark.}
What we considered in Example \ref{ex-sym} is 
a different kind of 
positive structure on $B^-_{\til w}$ where
$\til w=s_1s_2\cd s_n$.
More precisely, we define a rational morphism:
\[
\begin{array}{cccc}
\til \theta: &(\bbC^\tm)^n & 
\longrightarrow & B^-_{\til w}\\
&(a_1,\cd,a_n)& \mapsto & 
y_1(\frac{1}{c_1})\al_1^\vee(c_1)\cd 
y_n(\frac{1}{c_n})\al_n^\vee(c_n),
\end{array}
\]
where $c_i=a_1a_2\cd a_i$. Then it is easy to 
see that $\til\theta$ gives a positive structure
on $B^-_{\til w}$. Indeed, 
since we have
\[
e^c_{i}(Y_{\til w}(c_1,\cd,c_n))
=Y_{\til w}(c_1,\cd,c_{i-1,}cc_i,c_{i+1},\cd,c_n),
\]
we obtain 
\[
e_{i,\til\theta}(c,(a_1,\cd,a_n,a_{n+1}))
=
(a_1,\cd,ca_i,c^{-1}a_{i+1},\cd,a_n,a_{n+1}),
\]
where $a_1\cd a_{n+1}=1$.
The ultra-discritization of the geometric crystal 
on $B^-_{\til w}$ with respect to $\til\theta$ is as 
follows;
Set $\til B:=
\{(x_1,\cd,x_{n+1})\in \ZZ^{n+1}|x_1+\cd+x_{n+1}=0\}$ 
and for $x:=(x_1,\cd,x_{n+1})\in \til B$,
set
\[
\til e_i^c(x)=
(x_1,\cd,x_i+c,x_{i+1}-c,\cd,x_{n+1}) \q
(c\in \ZZ_{\geq0}),
\]
and $\til f_i^c=\til e_i^{-c}$.
Then ${\cal UD}_{\til\theta,\bbC^n}
(B^-_{\til w},\gamma,\{e_i\})$ is the Langlands dual of 
the crystal $\til B$.
The crystal $\til B$ holds the similar structure to 
some limit of ``crystal base for the
symmetric tensor module''.

\renewcommand{\thesection}{\arabic{section}}
\section{Tropical Braid-type isomorphisms}
\setcounter{equation}{0}
\renewcommand{\theequation}{\thesection.\arabic{equation}}
\newtheorem{pro6}{Proposition}[section]
\newtheorem{thm6}[pro6]{Theorem}
\newtheorem{ex6}[pro6]{Example}
\newtheorem{lem6}[pro6]{Lemma}

As an application of the tropicalization/ultra-discritzation
given in the previous section, 
we shall give a new proof of 
 the braid-type isomorphisms of 
crystals \cite{N}.
In order to do it, let us give the 
''tropical braid-type isomorphism''
(similar formula is given in \cite{BZ}):
\begin{pro6}
\label{trop-braid}
We have the following identities:
\begin{enumerate}
\item Type $A_2$:
\begin{eqnarray}
&& \hspace{-40pt}y_i(\frac{1}{c_1})\al^\vee_i(c_1)
y_j(\frac{1}{c_2})\al^\vee_j(c_2)
y_i(\frac{1}{c_3})\al^\vee_i(c_3)
\label{tro-A2} \\
&& \hspace{-40pt} =
y_j(\frac{c_1}{c_1c_3+c_2})
\al^\vee_j(\frac{c_1c_3+c_2}{c_1})
y_i(\frac{1}{c_1c_3})\al^\vee_i(c_1c_3)
y_j(\frac{c_1c_3+c_2}{c_1c_2})
\al^\vee_j(\frac{c_1c_2}{c_1c_3+c_2})
\nn
\end{eqnarray}
\item Type $B_2$ \q
$(\lan \al_i^\vee,\al_j\ran=-2,\,
\lan \al_j^\vee,\al_i\ran=-1)$:
\begin{eqnarray*}
&& y_i(\frac{1}{c_1})\al^\vee_i(c_1)
y_j(\frac{1}{c_2})\al^\vee_j(c_2)
y_i(\frac{1}{c_3})\al^\vee_i(c_3)
y_j(\frac{1}{c_4})\al^\vee_j(c_4)
\\
&& \qq\q=
y_j(\frac{1}{d_1})\al^\vee_j(d_1)
y_i(\frac{1}{d_2})\al^\vee_i(d_2)
y_j(\frac{1}{d_3})\al^\vee_j(d_3)
y_i(\frac{1}{d_4})\al^\vee_i(d_4),
\end{eqnarray*}
\begin{eqnarray}
{\rm where\,\, }
& d_1=c_4
+\frac{1}{c_2}\left(c_3+\frac{c_2}{c_1}\right)^2,\q
& d_2=c_1c_4+c_3+\frac{c_1c_3^2}{c_2},
\label{tro-B2-1}\\
& d_3=\frac{1}
{\frac{1}{c_2}+
\frac{1}{c_2^2c_4}\left(c_3+\frac{c_2}{c_1}\right)^2},\q
& d_4=\frac{1}
{\frac{c_4}{c_3}+\frac{c_3}{c_2}+\frac{1}{c_1}}.
\label{tro-B2-2}
\end{eqnarray}
\item Type $G_2$ \q
$(\lan \al_i^\vee,\al_j\ran=-3,\,
\lan \al_j^\vee,\al_i\ran=-1):$
\begin{eqnarray}
&& \hspace{-45pt}y_i(\frac{1}{c_1})\al^\vee_i(c_1)
y_j(\frac{1}{c_2})\al^\vee_j(c_2)
y_i(\frac{1}{c_3})\al^\vee_i(c_3)
y_j(\frac{1}{c_4})\al^\vee_j(c_4)
y_i(\frac{1}{c_5})\al^\vee_i(c_5)
y_j(\frac{1}{c_6})\al^\vee_j(c_6)\nn
\\
&& \hspace{-45pt}=
y_j(\frac{1}{d_1})\al^\vee_j(d_1)
y_i(\frac{1}{d_2})\al^\vee_i(d_2)
y_j(\frac{1}{d_3})\al^\vee_j(d_3)
y_i(\frac{1}{d_4})\al^\vee_i(d_4)
y_j(\frac{1}{d_5})\al^\vee_j(d_5)
y_i(\frac{1}{d_6})\al^\vee_i(d_6),\nn \\
&&  \label{G-rel}
\end{eqnarray}
where
\begin{eqnarray}
&& \hspace{-50pt}d_1 = 
\frac{1}{c_2^2}\left({c_3}+\frac{c_2}{c_1}\right)^3
+\frac{1}{c_4}
\left({c_5}+\frac{c_4}{c_3}\right)^3
+ \frac{2 c_4}{c_2} +\frac{3 c_4}{c_1 c_3} + 
\frac{3 c_5}{c_1} 
+\frac{3 c_3 c_5}{c_2} 
+ c_6,
\label{d1}\\
&& \hspace{-50pt}d_2 = 
\frac{c_1}{c_4}\left({c_5}+\frac{c_4}{c_3}\right)^3
+\frac{c_1c_3}{c_2^3}\left({c_3}+\frac{c_2}{c_1}\right)^3
+\frac{3c_1c_3c_5}{c_2}+\frac{2c_1c_4}{c_2}+
\frac{2c_4}{c_3}+c_1c_6+2c_5,
\nn \\
&& \label{d2}\\
&& \hspace{-50pt}d_5=
\frac{1}
{\frac{1}{c_6}\left(\frac{1}{c_4}\left(
{c_5}+\frac{c_4}{c_3}\right)^2+
\frac{c_3}{c_2}+\frac{1}{c_1}\right)^3+
\frac{c_6}{c_4}+\frac{3c_3c_5}{c_2c_4}
+\frac{3c_5}{c_1c_4}+\frac{3}{c_1c_3}
+\frac{2}{c_2}},\label{d5}\\
&& \hspace{-50pt}d_6=
\frac{1}
{\frac{1}{c_1}+\frac{c_3}{c_2}+
\frac{1}{c_4}\left(c_5+\frac{c_4}{c_3}\right)^2
+\frac{c_6}{c_5}},\qq
 d_3=\frac{c_2c_4c_6}{d_1d_5},\qq
 d_4=\frac{c_1c_3c_5}{d_2d_6}\label{d6}
\end{eqnarray}
\end{enumerate}
\end{pro6}

{\sl Proof.}
To show the above proposition, we need the following
well-known facts.
\begin{lem6} We have the following identities:
\begin{enumerate}
\item Type $A_2:$
\begin{equation}
y_i(a)y_j(b)=y_{\al_i+\al_j}(ab)y_j(b)y_i(a).
\label{y-A2}
\end{equation}
\item Type $B_2$ \q
$(\lan \al_i^\vee,\al_j\ran=-2,\,
\lan \al_j^\vee,\al_i\ran=-1):$
\begin{eqnarray}
&& y_i(a)y_j(b) =  y_{2\al_i+\al_j}(a^2b)
y_{\al_i+\al_j}(ab)y_j(b)y_i(a),\\
&& y_i(a)y_{\al_i+\al_j}(b)
 = y_{2\al_i+\al_j}(2ab)y_{\al_i+\al_j}(b)y_i(a).
\end{eqnarray}
\item Type $G_2$ \q
$(\lan \al_i^\vee,\al_j\ran=-3,\,
\lan \al_j^\vee,\al_i\ran=-1):$
\begin{eqnarray}
&&\hspace{-50pt}
y_i(a)y_j(b) =  y_{3\al_i+2\al_j}(a^3b^2)
y_{3\al_i+\al_j}(a^3b)y_{2\al_i+\al_j}(a^2b)
y_{\al_i+\al_j}(ab)y_j(b)y_i(a),\qq
\label{g2-1}\\
&& \hspace{-50pt}
y_{\al_i+\al_j}(a)y_{2\al_i+\al_j}(b)
 = y_{3\al_i+2\al_j}(3ab)y_{2\al_i+\al_j}(b)
y_{\al_i+\al_j}(a),\label{g2-2}\\
&& \hspace{-50pt}
y_j(a)y_{3\al_i+\al_j}(b) =
y_{3\al_i+2\al_j}(-ab)
y_{3\al_i+\al_j}(b)y_j(a).\label{g2-3}
\end{eqnarray}
\end{enumerate}
\end{lem6}
By using these relations, immediately we obtain  
the $A_2$ and $B_2$ cases.
The $G_2$ case is quite complicated to obtain the 
explicit form of $d_j$'s. Using (\ref{ae}),(\ref{g2-1}),
(\ref{g2-2}) and (\ref{g2-3}), we can
write the both sides of 
(\ref{G-rel}) in the form:
\[
y_{3\al_i+2\al_j}(A)
y_{3\al_i+\al_j}(B)y_{2\al_i+\al_j}(C)
y_{\al_i+\al_j}(D)y_j(E)y_i(F)
\al^\vee_i(G) \al^\vee_j(H).
\]
Then comparing the both sides, we get 
(\ref{d1}), (\ref{d2}), (\ref{d5}) and (\ref{d6}).\qed

By Proposition \ref{trop-braid}, we easily see that
each $d_j$ is a positive rational function in $c_j$'s.
Thus, the map 
\[
(c_1,c_2,\cd)\mapsto y_j(\frac{1}{d_1})\al_j^\vee(d_1)
y_i(\frac{1}{d_2})\al_i^\vee(d_2)\cd
\]
gives rise to positive structures on 
$B^-_{w_0}$ where $w_0$ is the longest element of the Weyl group of type
$A_2$, $B_2$ or $G_2$.
Then, if we consider the ultra-discritization of 
this positive strucutres, we obtain the so-called 
''braid-type isomorphisms'' between the tensor products
 of the crystal $B_i$'s(\cite{N}).

\begin{pro6}[\cite{N}]
\label{braid}
\begin{enumerate}
\item
If $\lan \al_i^\vee,\al_j\ran=
\lan \al_j^\vee,\al_i\ran=0$, 
$$
\begin{array}{c}
\phi^{(0)}_{ij}:
 B_i\ot B_i  \mapright{\sim}  B_j\ot B_i\\
\qq\q (x)_i\ot(y)_j  \mapsto (y)_j\ot (x)_i.
\end{array}
$$
\item
If $\lan \al_i^\vee,\al_j\ran=
\lan \al_j^\vee,\al_i\ran=-1$, 
$$
\phi^{(1)}_{ij}
:B_i\ot B_j\ot B_i\mapright{\sim} B_j\ot B_i\ot B_j,
$$
\begin{equation}
\hspace{-10pt}(z_1)_i\ot(z_2)_j\ot(z_3)_i 
\mapsto ({\rm max}(z_3,z_2-z_1))_j\ot (z_1+z_3)_i
\ot (-{\rm max}(-z_1,z_3-z_2))_j,
\label{A2-braid}
\end{equation}
\item
If $\lan \al_i^\vee,\al_j\ran=-1$, 
$\lan \al_j^\vee,\al_i\ran=-2$, 
$$
\phi^{(2)}_{ij}
:B_i\ot B_j\ot B_i\ot B_j
\mapright{\sim} B_j\ot B_i\ot B_j\ot B_i,
$$
$$
 (z_1)_i\ot(z_2)_j\ot(z_3)_i\ot (z_4)_j \mapsto
 (Z_1)_j\ot (Z_2)_i\ot (Z_3)_j\ot (z_4)_i
$$
\begin{equation}
\hspace{-30pt}\left\{
\begin{array}{l}
Z_1= \max(z_4,z_2-2z_1,2z_3-z_2)
\\
Z_2  = \max(z_1+z_4,z_3,z_1-z_2+2z_3)
\\
Z_3  =-\max(-z_2,-z_4-2z_1,-2z_2+2z_3-z_4)
\\
Z_4 =-\max(-z_3+z_4,-z_1,z_3-z_2)
\end{array}\right.
\label{B2-braid}
\end{equation}

\item
If $\lan \al_i^\vee,\al_j\ran=-1$, 
$\lan \al_j^\vee,\al_i\ran=-3$, 
$$
\hspace{-46pt}\phi^{(3)}_{ij}
:B_i\ot B_j\ot B_i\ot B_j\ot B_i\ot B_j
\mapright{\sim} 
B_j\ot B_i\ot B_j\ot B_i\ot B_j\ot B_i
$$
\beqnn
&&\hspace{-30pt}(z_1)_i\ot (z_2)_j\ot (z_3)_i\ot 
(z_4)_j\ot (z_5)_i\ot (z_6)_j\\
&& \qq \mapsto (Z_1)_j\ot (Z_2)_i\ot (Z_3)_j\ot 
(Z_4)_i\ot (Z_5)_j\ot (Z_6)_i,
\eeqnn
\begin{equation}
\hspace{-30pt}\left\{
\begin{array}{l}
Z_1=\max(z_6,3z_5-z_4,-3z_3+2z_4,-2z_2+3z_3,-3z_1+z_2)
\\
Z_2=\max(z_1+z_6,z_1-z_4+3z_5,
z_1-3z_3+2z_4,z_1-2z_2+3z_3,-z_1+z_3)
\\
Z_3  =  z_2+z_4+z_6-Z_1-Z_5,\\
Z_4 =  z_1+z_3+z_5-Z_2-Z_6,\\
Z_5=-\max(-z_4+z_6,-3z_4+6z_5-z_6,-6z_3+3z_4-z_6,\\
\qq\qq\qq\qq\qq \qq -3z_2+3z_3-z_6,-3z_1-z_6)
\\
Z_6=-\max(-z_1,-z_2+z_3,-z_4+2z_5,-2z_3+z_4,-z_5+z_6)
\end{array}\right.
\label{G2-braid}
\end{equation}
\end{enumerate}
\end{pro6}
We call $\phi^{(k)}_{ij}$ $(k=0,1,2,3)$ a {\it braid-type isomorphism}.

{\sl Proof.}
The formula in (\ref{A2-braid}),
(\ref{B2-braid}) and (\ref{G2-braid}) are obtained by rewriting 
the ones in \cite[Proposition 4.1]{N}
by using:
\begin{eqnarray*}
&\hspace{-90pt}
a_1+(a_2+(a_3+(\cd+(a_k)_+\cd )_+)_+)_+\\
&\qq\qq
={\rm max}(a_1,a_1+a_2,a_1+a_2+a_3,\cd,a_1+\cd +a_k).
\end{eqnarray*}
In (\ref{tro-A2}), the ultra-discritizations of 
$\frac{c_1c_3+c_2}{c_1}$, 
$c_1c_3$ and 
$\frac{c_1c_2}{c_1c_3+c_2}$ are 
$$
\begin{array}{lll}
v(\frac{c_1c_3+c_2}{c_1}) & = &
\max(v(c_1)+v(c_3),v(c_2))-v(c_1)
=\max(v(c_3),v(c_2)-v(c_1)),\\
v(c_1c_3) & =& v(c_1)+v(c_3),\\
v(\frac{c_1c_2}{c_1c_3+c_2})
&=& v(c_1)+v(c_2)-\max(v(c_1)+v(c_3),v(c_2))
=-\max(v(c_3)-v(c_2),-v(c_1))
\end{array}
$$
Thus, replacing $v(c_i)$ with $z_i$, we obtain 
(\ref{A2-braid}).

Similarly, considering 
the ultra-discritizations of $d_i$'s
in (\ref{tro-B2-1}) and (\ref{tro-B2-2}), we get 
(\ref{B2-braid}). Here note that 
in Proposition \ref{braid} (iii), we suppose 
$\lan \al_i^\vee,\al_j\ran=-1$, 
$\lan \al_j^\vee,\al_i\ran=-2$, which is the
Langlands dual of the condition in 
Proposition \ref{trop-braid} (ii).

In order to get the formula (\ref{G2-braid}),
we consider, {\it e.g.}, $v(d_1):$
\begin{equation}
\begin{array}{lll}
v(d_1) &= & 
\max(-2z_2+3z_3,-3z_1+z_2,3z_5-z_4,-3z_3+2z_4,z_6,\\
& &\qq z_4-z_2,z_4-z_1-z_3,z_5-z_1,z_3+z_5-z_2)\,\q
(v(c_j)=z_j),
\end{array}
\label{UD-d1}
\end{equation}
which seems to be different from $Z_1$ in (\ref{G2-braid}).
But, it is easy to see 
that both are same by the following simple formula:

For $m_1,\cd,m_k\in \bbR$ and $t_1,\cd,t_k\in \bbR_{\geq0}$
satisfying $t_1+\cd+t_k=1$, we have
\[
\max\left(m_1,\cd,m_k,\sum_{j=1}^k t_jm_j\right)
=\max(m_1,\cd,m_k) 
\]
Indeed, in (\ref{UD-d1}) we have
\[
 \begin{array}{lll}
&&z_4-z_2= \frac{1}{2}A_1+\frac{1}{2}A_4, \qq
z_4-z_1-z_3 = \frac{1}{6}A_1+\frac{1}{3}A_2+
\frac{1}{2}A_4, \\
&&z_5-z_1 =  \frac{1}{6}A_1+\frac{1}{3}A_2
+\frac{1}{3}A_3+\frac{1}{6}A_4,\qq
z_3+z_5-z_2 =
\frac{1}{2}A_1+\frac{1}{3}A_3+\frac{1}{6}A_4,
 \end{array}
\]
where $A_1:=-2z_2+3z_3$, 
$A_2:=-3z_1+z_2$, $A_3:=3z_5-z_4$ and
$A_4:=-3z_3+2z_4$.

Hence we have $Z_1=v(d_1)$. Others are obtained
similarly. Thus, considering the Langlands dual,
we get the desired result.\qed

\end{document}